\documentclass{elsarticle}

\usepackage{hyperref}

\usepackage{amssymb,amsmath,amsfonts}

\newtheorem{theorem}{Theorem}[section]

\newtheorem{corollary}[theorem]{Corollary}
\newtheorem{lemma}[theorem]{Lemma}

\newtheorem{definition}[theorem]{Definition}

\newtheorem{example}{Example}[section]
\newtheorem{remark}{Remark}[section]

\def\normtwo#1{\Vert#1\Vert_2}

\newcommand{\eps}{u}

\begin{document}

\begin{frontmatter}

\title{Numerical stability of the  symplectic $LL^T$ factorization}

\author[MB]{Maksymilian Bujok\corref{cor1}}
\ead{mbujok@swps.edu.pl}

\author[MR]{Miroslav Rozlo\v{z}n\'{\i}k\corref{cor2}}
\ead{miro@math.cas.cz}

\author[AgS]{Agata Smoktunowicz\corref{cor4}}
\ead{A.Smoktunowicz@ed.ac.uk}

\author[AS]{Alicja Smoktunowicz\corref{cor3}}
\ead{alicja.smoktunowicz@pw.edu.pl}

\cortext[cor1]{Principal corresponding author}

\address[MB]{SWPS University of Social Sciences and Humanities,  Warsaw, Poland}
\address[MR]{Institute of Mathematics, Czech Academy of Sciences, Prague, Czech Republic}
\address[AgS]{School of Mathematics, University of Edinburgh,  Edinburgh, Scotland, United Kingdom}
\address[AS]{Faculty of Mathematics and Information Science, Warsaw University of Technology, Warsaw, Poland}

\begin{abstract}
In this paper we give the detailed error analysis of  two algorithms $W_1$ and $W_2$  for computing  the symplectic factorization
of a symmetric positive definite and symplectic matrix  $A  \in \mathbb R^{2n \times 2n}$ in the form $A=LL^T$,
 where  $L  \in \mathbb R^{2n \times 2n}$ is a symplectic block  lower triangular matrix.
We prove that Algorithm $W_2$ is numerically stable for a broader class of symmetric positive definite matrices $A \in \mathbb R^{2n \times 2n}$.
It means that  Algorithm $W_2$ is producing the computed factors $\tilde L$ in floating-point arithmetic with  machine precision $\eps$
such that $\normtwo{A-\tilde L {\tilde L}^T}=  {\cal O}(\eps \normtwo{A})$.
On the other hand, Algorithm $W_1$ is unstable, in general, for symmetric positive definite  and  symplectic matrix $A$.
In this paper we  also give corresponding bounds for  Algorithm $W_1$ that are weaker. We show that the factorization error depends on
the condition number $\kappa_2(A_{11})$ of the principal submatrix $A_{11}$.

Bounds for the loss of symplecticity of the lower block triangular matrices  $L$  for both Algorithms $W_1$ and $W_2$ that hold in exact arithmetic
for a broader class of symmetric positive definite matrices  $A$ (but not necessarily symplectic) are also   given.
The tests performed in \textsl{MATLAB} illustrate that our error bounds for considered algorithms  are reasonably sharp.
\end{abstract}


\end{frontmatter}

\section{Introduction} 
We study numerical properties  of  two algorithms $W_1$ and $W_2$ from  \cite{Maks}.
We repeat  the relevant material  on Algorithms $W_1$ and $W_2$ without proofs, thus making our exposition self-contained.
All proofs and  \textsl{MATLAB} codes of these algorithms  can be  found in \cite{Maks}.
\begin{definition}\label{sympl}
$A \in \mathbb R^{2n \times 2n}$ is symplectic if it satisfies   $A^TJA=J$, where
\[
J=\left(
\begin{array}{cc}
 0  &   I \\
 -I &   0
\end{array}
\right),
\]
where $I \in \mathbb R^{n \times n} $ denotes the identity matrix.
\end{definition}

Our task is to compute  the  symplectic factorization $A=LL^T$  of a given  symmetric positive definite  and symplectic matrix $A \in  \mathbb R^{2n \times 2n}$, where
the matrix $L \in \mathbb R^{2n \times 2n}$ is symplectic and has a lower triangular block form  partitioned  conformally with $J$.
\begin{lemma}\label{lemacik1}
A  nonsingular block lower triangular  matrix $L\in \mathbb R^{2n \times 2n}$ of the form
\begin{equation}\label{L}
L=\left(
\begin{array}{cc}
 L_{11} &   0 \\
 L_{21} &  L_{22}
\end{array}
\right)
\end{equation}
is symplectic  if and only if  $L_{22}=L_{11}^{-T}$ and $L_{21}^T L_{11}=L_{11}^T L_{21}$. It follows that $L^T$ is then also symplectic.
\end{lemma}

Symplectic matrices form a Lie group under matrix multiplications. The product  of two symplectic matrices   is also a symplectic matrix.
The symplectic group is closed under transposition. If $A$ is symplectic then the inverse of $A$ equals  $A^{-1}=J^TA^TJ$, and it is also symplectic.

\medskip
Several efficient numerical methods for computing such specific factorization of symplectic matrices exist and are used in practical computations.
Symplectic matrices  arise in many applications, e.g.  in  various aspects of mathematical physics, including the application of symplectic block matrices to special relativity and
optimal control theory. The symplectic  $LL^T$  factorization can be used to compute the symplectic $QR$ factorization and the  Iwasawa decomposition of symplectic matrix.
For  more details we refer the reader to  \cite{Benzi}, \cite{Angelika}- \cite{Dopico}, \cite{Tam}-\cite{ Volker}.

In this paper, we first consider an implementation of the $HH^T$ factorization from \cite{Dopico} denoted as Algorithm $W_1$.
\begin{theorem} [{\cite[Theorem 5.2 and Corollary 2.3]{Dopico}}]
\label{HHT}
 Let  $A \in \mathbb R^{2n \times 2n}$  be a symmetric positive definite and symplectic matrix  partitioned  as
\begin{equation}\label{blockA}
A=\left(
\begin{array}{cc}
 A_{11}  &  A_{12} \\
 A_{12}^T &   A_{22}
\end{array}
\right),
\end{equation}
where  $A_{11} \in \mathbb R^{n \times n}$  and  $A_{22} \in \mathbb R^{n \times n}$.
If $A_{11}=L_{11} L_{11}^T$ is the Cholesky decomposition of $A_{11}$, then $A=L_1L_1^T$, where
$L_1$ is symplectic and has a form
\begin{equation}\label{sympL}
L_1=\left(
\begin{array}{cc}
 L_{11} &   0 \\
 (L_{11}^{-1} A_{12})^T &  L_{11}^{-T}
\end{array}
\right).
\end{equation}
Moreover, we have
\begin{equation}\label{invA11}
S=A_{11}^{-1},
\end{equation}
where $S=A_{22}-A_{12}^TA_{11}^{-1}A_{12}$ denotes the Schur complement of $A_{11}$ in $A$.
\end{theorem}

\medskip
The second Algorithm $W_2$ was  proposed in \cite{Maks} and it  is  based on the following theorem.
\begin{theorem}
\label{AlgW2}
Let  $A \in \mathbb R^{2n \times 2n}$  be a symmetric positive definite  matrix
partitioned conformally with $J$ as in (\ref{blockA}).
If $A_{11}=L_{11} L_{11}^T$ is the Cholesky decomposition of $A_{11}$, then
there is a unique matrix $L_2 \in \mathbb R^{2n \times 2n}$ such that $A=L_2L_2^T$, where $L_2$ has a form
\begin{equation}\label{formL}
L_2=\left(
\begin{array}{cc}
 L_{11}  &  0 \\
(L_{11}^{-1} A_{12})^T    &   U_{22}
\end{array}
\right),
\end{equation}
where  $U_{22}$   can be taken as the Reverse Cholesky factor of the Schur complement
$S=A_{22}-A_{12}^TA_{11}^{-1}A_{12}$ of $A_{11}$ in $A$. Then $S=U_{22} U_{22}^T$, where $U_{22}$
is the upper triangular matrix  with positive diagonal entries.
\end{theorem}

\medskip
\begin{remark}\label{Remark0}
We recall that the  Reverse Cholesky factorization gives a unique decomposition  $A=UU^T$ of a  symmetric  positive definite matrix $A \in  \mathbb R^{n \times n}$,
where $U \in  \mathbb R^{n \times n}$ is upper triangular matrix with positive diagonal entries.
This is  the  Cholesky factorization of the symmetric positive definite matrix $P^T A P$,
where $P$ is the permutation matrix comprising the identity matrix with its column in reverse order.
If $P^T A P=L L^T$, where $L$ is lower triangular (with positive diagonal entries),
then $A=UU^T$, with $U=P L P^T$ being  upper triangular (with positive diagonal entries).
\end{remark}

\medskip
\begin{itemize}
\item[] {\bf Algorithm $W_1$}

Given a symmetric positive definite and symplectic matrix $A \in   \mathbb R^{2n \times 2n}$ this
algorithm computes the  symplectic  $LL^T$ factorization $A=L_1 L_1^T$, where $L_1$ is defined  in  (\ref{sympL}).
\begin{enumerate}
\item $A_{11}=L_{11}L_{11}^T$ (Cholesky).
\item $L_{11}  L_{21}^T=A_{12}$ (forward  substitution).
\item $L_{11} X=I$ (forward substitution).
\item $U_{22}=X^{T}$.
\item $L_{22}=U_{22}$.
\end{enumerate}
\noindent
{\bf Cost:} $\frac{5}{3} n^3$ flops.

\item[] {\bf Algorithm $W_2$}

Given a symmetric positive definite and  symplectic matrix $A \in  \mathbb R^{2n \times 2n}$ this
algorithm computes the  symplectic  $LL^T$ factorization $A=L_2 L_2^T$, where $L_2$ is given in  (\ref{formL}).
\begin{enumerate}
\item $A_{11}=L_{11}L_{11}^T$ (Cholesky).
\item $L_{11}  L_{21}^T=A_{12}$ (forward  substitution).
\item $S=A_{22}- L_{21} L_{21}^T$ (Schur complement).
\item $S=U_{22}\, U_{22}^T$ (Reverse Cholesky).
\item $L_{22}=U_{22}$.
\end{enumerate}
\noindent
{\bf Cost:} $\frac{8}{3} n^3$ flops.
\end{itemize}

\medskip

What insights can we derive from applying Algorithms $W_1$ and $W_2$ to a  symmetric positive definite matrix
$A \in \mathbb R^{2n \times 2n}$  but  not  symplectic?
We see  that the blocks $L_{11}$ and $L_{21}$   are the same for both these algorithms
but the block $L_{22}=U_{22}$ can be different.  Clearly,  the matrix $L_2$ defined in  (\ref{formL})  is not symplectic (since $A=L_2 L_2^T$).
Notice that in this case the matrix $L_1$, defined in  (\ref{sympL}),  may not be symplectic. We recall that from Lemma \ref{lemacik1} it follows that $L_1$ is symplectic
if and only if $L_{21}^T L_{11}=L_{11}^T L_{21}$. 
We  illustrate it on the following example.
\begin{example}\label{przyklad1}
Let
\[A=\left(
\begin{array}{cc}
 A_{11}  &  A_{12} \\
 A_{12}^T &   A_{22}
\end{array}
\right)
=\left(
\begin{array}{cc|cc}
1 & 1 & 1 & 1\\
1 & 2 & 2 & 2\\
\hline
1 & 2 & 3 & 3\\
1 & 2 & 3 & 4
\end{array}
\right).
\]
It is easy to check that $A$ is a symmetric positive definite matrix but it is not symplectic.
Let  $L_0$ be the Cholesky factor of $A$. Then $A=L_0 L_0^T$, where
\[
L_0=\left(
\begin{array}{cc|cc}
1 & 0 & 0 & 0\\
1 & 1 & 0 & 0\\
\hline
1 & 1 & 1 & 0\\
1 & 1 & 1 & 1
\end{array}
\right).
\]
If $L_1$ is computed by Algorithm $W_1$ then we have $A +\Delta A = L_1 L_1^T$, where
\[
L_1=\left(
\begin{array}{cc|cc}
1 & 0 & 0 & 0\\
1 & 1 & 0 & 0\\
\hline
1 & 1 & 1 & -1\\
1 & 1 & 0 & 1
\end{array}
\right),
\quad
\Delta A =  \left(
\begin{array}{cc|cc}
0 & 0 & 0 & 0\\
0 & 0 & 0 & 0\\
\hline
0 & 0 & 1 & -2\\
0 & 0 & -2 & -1
\end{array}
\right).
\]

Moreover, we have
\[
L_1^T J L_1- J=
\left(
\begin{array}{cc|cc}
0 & 1 & 0 & 0\\
-1 & 0 & 0 & 0\\
\hline
0 & 0 & 0 & 0\\
0 & 0 & 0 & 0
\end{array}
\right),
\]
and
\[
(A+\Delta A)^TJ (A+\Delta A)-J =\left(
\begin{array}{cc|cc}
0 & 1 & 1 & 1\\
-1 & 0 & 0 & 0\\
\hline
-1 & 0 & 0 & 1\\
-1 & 0 & -1 & 0
\end{array}
\right).
\]

We see  that neither matrix $L_1$ nor $A+\Delta A= L_1  L_1^T$ is symplectic.

If $L_2$ is obtained by Algorithm $W_2$ then $A=L_2 L_2^T$, in which
\[
L_2=\left(
\begin{array}{cc|cc}
1 & 0 & 0 & 0\\
1 & 1 & 0 & 0\\
\hline
1 & 1 &  \frac{\sqrt{2}}{2} & \frac{\sqrt{2}}{2}\\
1 & 1 & 0 & \sqrt{2}
\end{array}
\right).
\]
If we apply the  Reverse Cholesky factorization then $A=U_0U_0^T$, where
\[
U_0=\left(
\begin{array}{cc|cc}
\frac{\sqrt{2}}{2} & \frac{\sqrt{6}}{6}  &  \frac{\sqrt{3}}{6} & \frac{1}{2}\\
0 & \frac{\sqrt{6}}{3} & \frac{\sqrt{3}}{3} & 1\\
\hline
0 & 0 &  \frac{\sqrt{3}}{2} & \frac{3}{2}\\
0 & 0 & 0 & 2
\end{array}
\right).
\]
\end{example}

The paper is organized as follows.
In this paper,  for any matrix $A \in \mathbb R^{n \times n}$,   $|A|= (\vert a_{ij} \vert)$ denotes the absolute
value of $A = (a_{ij})$,  $\normtwo{A}$ is  the spectral norm of $A$, and $\kappa_2(A)=\normtwo{A^{-1}} \cdot \normtwo{A}$ is the condition number of $A$.
Inequalities between matrices are understood to hold componentwise, i.e.\   for given $A, B \in \mathbb R^{n \times n}$   we write $|A| \leq |B|$
if  $|a_{ij}|  \leq |b_{ij}|$ for  $i, j=1, \ldots, n$.
We recall  that    $|A| \leq |B|$ implies $\normtwo{|A|} \leq \normtwo{|B|}$. We recall that  for any matrix $A \in \mathbb R^{n \times n}$ we have the bounds
\[
\normtwo{A} \leq \normtwo{|A|} \leq \sqrt{n} \, \normtwo{A}.
\]

In Section $2$,  we study some  properties of the symplectic $LL^T$ factorization
that hold in exact arithmetic.
In Section $3$, we consider finite precision arithmetic and  derive  bounds for the factorization error in  Algorithms $W_1$ and $W_2$.
We prove  that Algorithm $W_2$ is numerically stable for a broader class of symmetric positive definite matrices $A \in \mathbb R^{2n \times 2n}$.
It means that  the matrix $\tilde L_2$, computed  by Algorithm $W_2$  in floating-point arithmetic with  machine precision $\eps$, is the exact factor of the $LL^T$
factorization of o slightly perturbed matrix $A+ E$, where  $\normtwo{E}=  {\cal O}(\eps \normtwo{A})$, see Theorem \ref{ThmW2}.
Then  the factorization error satisfies $\normtwo{A-\tilde L_2 {\tilde L_2}^T}=  {\cal O}(\eps \normtwo{A})$.
We show that Algorithm $W_1$ is unstable, in general, for symmetric positive definite  and  symplectic matrix $A$. In this case,  the factorization error $\normtwo{A-\tilde L_1 {\tilde L_1}^T}$ depends on
the condition number $\kappa(A_{11})=\normtwo{A_{11}^{-1}} \cdot \normtwo{A_{11}}$ of the principal submatrix $A_{11}$, see Theorem \ref{NewThmW1} and Remark \ref{RemarkB}.

In Section $4$,  we present  numerical tests and comparison of the methods.
Concluding remarks  are given in Section $5$.

\medskip
\section{Theoretical properties of the symplectic $LL^T$ factorization}
In practice, we often cannot  ensure that the matrix  $A \in \mathbb R^{2n \times 2n}$  is  symplectic.
What follows from Theorem \ref{HHT} and Theorem \ref{AlgW2} when  Algorithms $W_1$  and $W_2$ are
applied to a broader class of symmetric positive definite matrices $A$ but not necessarily symplectic?

\begin{lemma}\label{lemacik}
Let   $A  \in \mathbb R^{2n \times 2n}$  be  a symmetric positive definite matrix (not necessarily  symplectic).
Assume that  $A_{11}=L_{11} L_{11}^T$ is the Cholesky decomposition of $A_{11}$, $L_1 \in \mathbb R^{2n \times 2n}$ is defined in  (\ref{sympL}), and
let  $A=L_{2}L_{2}^T$, where $L_2$ is given in  (\ref{formL}).
Then
\begin{equation}\label{deltaA}
A+ \Delta A = L_1 L_1^T, \quad  \Delta A=
\left(
\begin{array}{cc}
0 &  0 \\
0 & A_{11}^{-1}-S
\end{array}
\right),
\end{equation}
where  $S=A_{22}-A_{12}^TA_{11}^{-1}A_{12}$ denotes the Schur complement of $A_{11}$ in $A$.
Equivalently, we can formulate (\ref{deltaA}) in terms of the following block $LDL^T$ decompositions:
\begin{equation}\label{blockLDL1}
L_1 L_1^T =
\left(
\begin{array}{cc}
I &  0 \\
W^T & I
\end{array}
\right)
\left(
\begin{array}{cc}
A_{11} &  0 \\
0 &  A_{11}^{-1}
\end{array}
\right)
\left(
\begin{array}{cc}
I &  W \\
0 & I
\end{array}
\right),
\end{equation}
\begin{equation}\label{blockLDL2}
A=L_2 L_2^T =
\left(
\begin{array}{cc}
I &  0 \\
W^T & I
\end{array}
\right)
\left(
\begin{array}{cc}
A_{11} &  0 \\
0 & S
\end{array}
\right)
\left(
\begin{array}{cc}
I &  W \\
0 & I
\end{array}
\right),
\end{equation}
where the matrix $W$ is given as
$W=A_{11}^{-1} A_{12}$.
\end{lemma}
\noindent
\begin{proof}
Let $L_{21}= (L_{11}^{-1} A_{12})^T$ and  $L_{22}=L_{11}^{-T}$. Then  $L_{22} L_{22}^T=A_{11}^{-1}$.
Thus, $A_{22}-L_{21} L_{21}^T=S$ is  the Schur complement  of $A_{11}$ in $A$.
Then simple calculation shows
\[
\Delta A=
\left(
\begin{array}{cc}
L_{11}  L_{11}^T-A_{11} &  L_{11}  L_{21}^T-A_{12}  \\
 ( L_{11}  L_{21}^T- A_{12})^T & L_{21} L_{21}^T +  L_{22} L_{22}^T-A_{22}
\end{array}
\right)=
\left(
\begin{array}{cc}
0 &  0 \\
0 & A_{11}^{-1}-S
\end{array}
\right),
\]
which establishes the result (\ref{deltaA}).
Since $L_{21}=W^TL_{11}$, we have
\begin{equation}
\label{L1}
L_1 =\left(
\begin{array}{cc}
I &  0 \\
W^T & I
\end{array}
\right)
\left(
\begin{array}{cc}
L_{11} &  0 \\
0 & L_{11}^{-T}
\end{array}
\right), \quad
L_2=\left(
\begin{array}{cc}
I &  0 \\
W^T & I
\end{array}
\right)
\left(
\begin{array}{cc}
L_{11} &  0 \\
0 & U_{22}
\end{array}
\right ),
\end{equation}
where $S=U_{22}U_{22}^T$. This leads to  (\ref{blockLDL1})-(\ref{blockLDL2}).
\end{proof}

\medskip

\begin{remark}\label{UwagaA3}
{\bf Notice that immediately we have that ${\kappa_2(L_2)}^2= \kappa_2(A)$.
It is evident by Theorem \ref{HHT} that   $\Delta A=0$ and $L_1=L_2$  with $\kappa_2(L_1)= \kappa(L_2)$  for any symmetric positive definite and symplectic matrix $A$.
From   $A+\Delta A=L_1 L_1^T$ and (\ref{deltaA}) it follows that
\begin{equation}\label{Miro1}
\normtwo{L_1}^2 = \normtwo{A+\Delta A} \leq \normtwo{A} + \normtwo{\Delta A}= \normtwo{A}+\normtwo{A_{11}^{-1} - S}.
\end{equation}

To estimate $\normtwo{L_1^{-1}}^2$, we rewrite  $A+\Delta A$ as follows $A+\Delta A= L_2 (I+H)L_2^T$, where $H=L_2^{-1} \Delta A L_2^{-T}$.
Notice that $H$ is symmetric, hence  $\normtwo{H}=\rho(H)$, where $\rho(H)$ is  the spectral radius of  $H$.
We check immediately that
\[
\rho(H)=\rho(L_2^T (A^{-1} \Delta A) L_2^{-T})=\rho(A^{-1} \Delta A).
\]

Now assume that $\rho(A^{-1} \Delta A)<1$. Then $I+H$ is nonsingular and we get
\[
(A+\Delta A)^{-1}= L_1^{-T} L_1^{-1}= L_2^{-T} (I+H)^{-1}L_2^{-1}.
\]

This  leads to
\begin{equation}\label{Miro2}
\normtwo{L_1^{-1}}^2 \leq \normtwo{ L_2^{-1}}^2/(1-\normtwo{H}) = \normtwo{ A^{-1}}/(1-\rho(A^{-1} \Delta A)).
\end{equation}

Notice that
\[
A^{-1} \Delta A=
\left(
\begin{array}{cc}
0  &  -WS^{-1}(A_{11}^{-1}-S) \\
0  & S^{-1}(A_{11}^{-1}-S)
\end{array}
\right),
\]
which gives the formula  $\rho(A^{-1}\Delta A) = \rho(S^{-1}(A_{11}^{-1}-S))$.

This together with  $\kappa_2(A+\Delta A)={\kappa_2(L_1)}^2$ and  (\ref{Miro1})-(\ref{Miro2}) leads to 
\begin{equation}\label{Miro3}
\kappa_2(A+\Delta A) \leq \frac{\kappa_2(A)}{1-\rho(S^{-1}(A_{11}^{-1}-S))} \cdot  \left(1+\frac{\normtwo{A_{11}^{-1} - S }}{\normtwo{A}}\right).     
\end{equation}
}
\end{remark}

\medskip
\begin{remark}\label{UwagaA1}
Under the  hypotheses  of Lemma \ref{lemacik} we have
\begin{equation}\label{Important decW1}
\frac{\normtwo{A-L_1 L_1^T}}{\normtwo{A}}=\frac{\normtwo{\Delta A}}{\normtwo{A}}= \frac{\normtwo{A_{11}^{-1} - S }}{\normtwo{A}}.
\end{equation}
Notice that   $A_{11}^{-1}=L_{11}^{-T} L_{11}^{-1}$ and $S=A_{22}-L_{21} L_{21}^T$, where $L_{21}^T=L_{11}^{-1} A_{12}$ and $ L_{22}=L_{11}^{-T}$,   hence in practice
we can compute the factorization error  as follows
\begin{equation}\label{computeL1}
\frac{\normtwo{A-L_1 L_1^T}}{\normtwo{A}} =\frac{\normtwo{L_{22} L_{22}^T - (A_{22}-L_{21} L_{21}^T)}}{\normtwo{A}}.
\end{equation}
\end{remark}

\medskip

Now we show how to derive   the  perturbation results for the symplectic  $LL^T$ factorization  from those for the Cholesky factorization.
We apply here  norm invariant properties of orthogonal matrices. We recall that  $\normtwo{A} \leq ||A||_F \leq \sqrt{n} \normtwo{A}$
for any matrix $A \in \mathbb R^{n \times n}$,
where $||A||_F$ denotes the Frobenius norm of the matrix $A$.
If $Q  \in \mathbb R^{n \times n}$ is an orthogonal matrix then  $|| QA||_p= || AQ||_p=|| QAQ^T||_p=|| A|_p$ for $p=2, F$.

\begin{theorem}[{\cite[Theorem 1.4]{Sun}}]
\label{Sun}
Let $A \in \mathbb R^{n \times n}$ be symmetric positive definite with the Cholesky factorization $A=L  L^T$ and let $E$ be a symmetric matrix satisfying
$\normtwo{A^{-1}} \normtwo{E}<1$. Then $A+E$ has the  Cholesky factorization $A+E= (L+\Delta L) (L+\Delta L)^T$, where
 \begin{equation}\label{ThmSun}
 \frac{||\Delta L||_{F}}{|| L||_p} \leq  2^{-1/2}   \, \frac{\kappa_2(A)}{1-\normtwo{A^{-1}} \normtwo{E}} \cdot  \frac{||E||_{F}}{|| A||_p},  \quad p=2, F.
\end{equation}
\end{theorem}

\medskip
\noindent
Similar result can be obtained for the Reverse Cholesky factorization.
\noindent
\begin{theorem}\label{ReChol1}
Let $A \in \mathbb R^{n \times n}$ be symmetric positive definite with the  Reverse Cholesky factorization $A=U U^T$, where $U$ is upper triangular,   and let $E$ be a symmetric matrix satisfying
$\normtwo{A^{-1}} \normtwo{E}<1$. Then $A+E$ has the Reverse Cholesky factorization $A+E= (U+\Delta U) (U+\Delta U)^T$, where
\begin{equation}\label{U_Reverse}
\frac{||\Delta U||_{F}}{|| U||_p} \leq  2^{-1/2}   \, \frac{\kappa_2(A)}{1-\normtwo{A^{-1}} \normtwo{E}} \cdot  \frac{||E||_{F}}{|| A||_p},  \quad p=2, F.
\end{equation}
\end{theorem}
\noindent
\begin{proof} The Reverse Cholesky is, in fact,  the  Cholesky factorization of the symmetric positive definite matrix $P^T A P$,
where $P$ is the permutation matrix comprising the identity matrix with its column in reverse order (see Remark \ref{Remark0}). We recall that these factorizations are unique.

Let $A$ has the Reverse Cholesky factorization $A=UU^T$. Then the matrix  $PAP^T$ has the Cholesky factorization $PAP^T=L L^T$,
where $L= P^T U P$. From Theorem \ref{Sun} it follows that there exists a lower triangular matrix $\Delta L$ such that
$P(A+E)P^T=PAP^T+ PEP^T= (L+\Delta L)(L+\Delta L)^T$ with the bounds  (\ref{ThmSun}).
Since $P$ is orthogonal, then we have the identities $||A||_p=||PAP^T||_p$, $||A^{-1}||_p=||(PAP^T)^{-1}||_p$  and $||E||_p=||PEP^T||_p$ for $p=2, F$.
We conclude that $A+E$ has the Reverse Cholesky factorization $A+E=(U+\Delta U) (U+\Delta U)^T$, where $\Delta U=P\Delta L P^T$ is upper triangular.
Clearly, $||\Delta L||_p=||\Delta U||_p$, which establishes the formula (\ref{U_Reverse}).
The proof is now complete.
\end{proof}

\medskip
\noindent
The following  theorem  shows how the matrix $L_2$, defined in Theorem \ref{AlgW2}, changes when $A$ is perturbed.
\begin{theorem}\label{perturb_L2}
Let $A \in \mathbb R^{2n \times 2n}$ be symmetric positive definite and  with the  factorization $A=L_2 L_2^T$, where
$L_2$ is defined in (\ref{formL}), and let $E$ be a symmetric matrix satisfying
$\normtwo{A^{-1}} \normtwo{E}<1$. Then $A+E$ has the  factorization
$A+E= (L_2+\Delta L_2) (L_2+\Delta L_2)^T$, where $\Delta L_2$ has the same form as $L_2$, and satisfies
\begin{equation}\label{perturbL2}
\frac{||\Delta L_2||_{F}}{|| L_2||_p} \leq  2^{-1/2}   \, \frac{\kappa_2(A)}{1-\normtwo{A^{-1}} \normtwo{E}} \cdot  \frac{||E||_{F}}{|| A||_p},  \quad p=2, F.
\end{equation}
\end{theorem}
\noindent
\begin{proof}
Our proof is similar to that of  Theorem \ref{ReChol1}.
Let $Q=diag(I, P)$, where $P$  is the permutation matrix comprising the identity matrix $I$ with its column in reverse order.
Then the matrix $Q A Q^T$ is symmetric positive definite and with the  factorization $Q A Q^T=L L^T$, where
\[
L=Q  L_2 Q^T=
\left(
\begin{array}{cc}
  L_{11}  &  0 \\
PL_{21} &   PL_{22}P^T
\end{array}
\right).
\]

Notice that $PL_{22}P^T$ is lower triangular, hence $L$ is also a lower triangular matrix. Thus, we get the Cholesky factorization of the matrix $Q A Q^T$.
Now we can apply Theorem \ref{Sun}. It   follows that there exists a lower triangular matrix $\Delta L$ such that
$Q(A+E)Q^T=QAQ^T+ QEQ^T= (L+\Delta L)(L+\Delta L)^T$ with the bounds  (\ref{ThmSun}).
Then $A+E= (L_2+\Delta L_2) (L_2+\Delta L_2)^T$, where $\Delta L_2=Q \Delta L Q^T$ has the same form as the matrix $L_2$, which
proves  (\ref{perturbL2}).
\end{proof}

\medskip

{\bf We define the loss of symplecticity of  arbitrary matrix $A \in \mathbb R^{2n \times 2n}$ as  $\normtwo{\Omega(A)}$, where $\Omega(A) = A^TJA-J$.
In practice, sometimes small perturbations can result in a significant deviation of a given symplectic  matrix $A$ from symplecticity. 
For this reason, theoretical properties of algorithms may diverge from computational properties.
If $A$ is symplectic then $\Omega(A) =0$, $A^{-1}=-JA^TJ$, hence $\normtwo{A^{-1}}=\normtwo{A}$ and $\kappa_2(A)=\normtwo{A}^2$.
Let  $\hat A=A+E$ denote the perturbed symplectic matrix  $A \in \mathbb R^{2n \times 2n}$. Then  $\Omega(\hat A) =\Omega(A) + (A^TJE+E^TJA+E^TJE)$, hence we have 
$\normtwo{\Omega(\hat A)} \leq 2 \normtwo{A} \normtwo{E}+ \normtwo{E}^2$. If $\normtwo{E} \leq \epsilon \normtwo{A}$ with $\epsilon <1$ then 
$\normtwo{\Omega(\hat A)} \leq (2\epsilon +\epsilon^2) \normtwo{A}^2=  (2\epsilon +\epsilon^2) \kappa_2(A)$. We show in the following example that if the matrix $A$ is ill-conditioned, meaning the condition number $\kappa_2(A)$ is large, then the loss of $\hat A$ from symplecticity may be significant. 
}

\begin{example}\label{macierz2}
{\bf  Let
\[
A=\left(
	\begin{array}{cc}
	G  &  I \\
	I &   2 G^{-1}
	\end{array}
	\right),
	\quad
\hat A=\left(
	\begin{array}{cc}
	\hat G  &  I \\
	I &   2 G^{-1}
	\end{array}
	\right),	
\]	
where $G=diag(t,1/t)$ and  $\hat G=diag(t,1/t+\theta)$, where $t \ge 1$ and $\theta >0$. Clearly,  $A$ is  symmetrix positive definite and symplectic  (see also Example \ref{example3}).
We have  $\hat A=A+E$, where $E=diag(0,\theta,0,0)$. Clearly,  $\normtwo{E}=\theta$ and $\normtwo{\Omega(\hat A)}=2t \theta$. It is easy to check that 
$2 t  \leq \normtwo{A} \leq 2t+1$, hence $4 t^2  \leq \kappa_2(A) \leq (2t+1)^2$. We see that the matrix $A$ is  ill-conditioned for large values of $t$.
For example, if we take  $t={10}^6$ and $\theta=10^{-10}$ then $\kappa_2(A) \ge 4 t^2= 4 \cdot {10}^{12}$, and the loss $\hat A$ from  symplecticity equals $\normtwo{\Omega(\hat A)}=2 \cdot {10}^{-4}$.}
\end{example}

It is interesting to see how the quantity $\normtwo{\Delta A}=\normtwo{A_{11}^{-1} - S}$ from Lemma \ref{lemacik} depends on the other  quantity $\normtwo{\Omega(A)}$
(the loss of symplecticity).

\medskip
\begin{lemma}
\label{sympA}
Let  $A \in \mathbb R^{2n \times 2n}$  be a symmetric positive definite  matrix  partitioned conformally with $J$ as
in (\ref{blockA}).  Then
\begin{equation}\label{OmegaA}
\Omega(A) =\left(
\begin{array}{cc}
 \Omega_{11}  &  \Omega_{12} \\
\Omega_{21}   &   \Omega_{22}
\end{array}
\right), \quad \Omega_{21}=-\Omega_{12}^T,
\end{equation}
where
\begin{eqnarray}\label{Omega11}
\Omega_{11} & = & A_{11} A_{12}^T-A_{12} A_{11},
\\
\label{Omega12}
\Omega_{12} & = & A_{11} A_{22}-A_{12}^2-I,
\\
\label{Omega22}
\Omega_{22}& =  & A_{12}^T A_{22}-A_{22} A_{12}.
\end{eqnarray}

\medskip
Let $S$ be the Schur complement of $A_{11}$ in $A$ and $W=A_{11}^{-1} A_{12}$.
Then
\begin{equation}\label{M1}
\Omega_{11}= A_{11} (W^T-W) A_{11}
\end{equation}
and
\begin{equation}\label{M2}
A_{11}^{-1}-S=A_{11}^{-1} \left (\Omega_{11} W -\Omega_{12}\right).
\end{equation}
\end{lemma}
\noindent
\begin{proof}
From (\ref{Omega11}) it follows that $\Omega_{11}= A_{11} (A_{12}^T A_{11}^{-1}-A_{11}^{-1}A_{12}) A_{11}$, which completes the proof of (\ref{M1}).
To prove (\ref{M2}),  we rewrite $A_{11}^{-1}-S$ as follows $A_{11}^{-1}-S=A_{11}^{-1} (I-A_{11}S)$.
Since $S=A_{22}-A_{12}^T A_{11}^{-1} A_{12}$,  we have
\begin{equation}\label{M3}
A_{11}S=A_{11} A_{22}-(A_{11}A_{12}^T) A_{11}^{-1} A_{12}.
\end{equation}
From  (\ref{Omega12}) it follows that
\begin{equation}\label{M4}
A_{11} A_{22}=I+A_{12}^2+ \Omega_{12}.
\end{equation}
Next, by  (\ref{Omega11}) we get
\begin{equation}\label{M5}
A_{11}A_{12}^T=A_{12}A_{11}+ \Omega_{11}.
\end{equation}
Substituting  (\ref{M4})-(\ref{M5}) to (\ref{M3}) leads to  $A_{11}S=I+\Omega_{12}- \Omega_{11} (A_{11}^{-1}A_{12})$.
Thus, $A_{11}^{-1}-S=A_{11}^{-1} (I-A_{11}S)=A_{11}^{-1} \left (\Omega_{11} W -\Omega_{12}\right)$. The proof is complete.
\end{proof}

\begin{remark}\label{UwagaA2}
Under the  hypotheses  of Lemma \ref{sympA} we get the bound
\[
 \normtwo{A_{11}^{-1}-S} \leq  \normtwo{A_{11}^{-1}} \left(\normtwo{\Omega_{11}(A)} \normtwo{W}+ \normtwo{\Omega_{12}(A)}\right).
\]
Due to $\normtwo{\Omega_{11}} \leq  \normtwo{\Omega(A)}$ and $\normtwo{\Omega_{12}} \leq  \normtwo{\Omega(A)}$
we obtain the inequality
\begin{equation}\label{T2}
\normtwo{A_{11}^{-1}-S} \leq \normtwo{A_{11}^{-1}} \left(\normtwo{W}+ 1 \right)   \normtwo{\Omega(A)}.
\end{equation}
Indeed, the quantity $\normtwo{\Delta A}= \normtwo{A_{11}^{-1}-S} $, defined in Lemma \ref{lemacik}, is bounded up to the factor
$\normtwo{A_{11}^{-1}} \left(\normtwo{W}+ 1 \right)$ by the quantity $\normtwo{\Omega(A)}$.
This essentially shows that  in fact $\normtwo{\Delta A}$ in some way also measures
a departure of the matrix $A$ from   symplecticity.
\end{remark}

\medskip
\noindent
Note that the matrix $W=A_{11}^{-1} A_{12}$ is present  in the results stated in  Lemmas \ref{lemacik} and  \ref{sympA}.
The following lemma demonstrates how to bound $\normtwo{W}$.

\begin{lemma} [{\cite[Lemma 10.12]{higham:2002}}]
\label{boundW}
If $A  \in \mathbb R^{2n \times 2n}$, partitioned as in  (\ref{blockA}),  is  symmetric positive definite then  $W=A_{11}^{-1} A_{12}$  satisfies
\[
\normtwo{W}^2 \leq \normtwo{A_{11}^{-1}} \normtwo{A_{22}}.
\]
\end{lemma}

\medskip
\noindent
Next lemma presents  the perturbation bounds for the Schur complement.
\begin{lemma}[{\cite[Lemma 10.10]{higham:2002}}]
\label{bound_Schur}
Let $A  \in \mathbb R^{2n \times 2n}$  be symmetric and positive definite  partitioned as in  (\ref{blockA}).
Let $\tilde A= A+\tilde E$, where $\tilde E  \in \mathbb R^{2n \times 2n}$  is symmetric.
Assume that  $\tilde A$ is positive definite and  is  partitioned  conformally with $A$ as
\[
\tilde A=\left(
\begin{array}{cc}
\tilde  A_{11}  & \tilde  A_{12} \\
\tilde  A_{12}^T &   \tilde A_{22}
\end{array}
\right),
\quad
\tilde E=\left(
\begin{array}{cc}
\tilde  E_{11}  & \tilde  E_{12} \\
\tilde  E_{12}^T &   \tilde E_{22}
\end{array}
\right).
\]
Let $S=A_{22}-A_{12}^T A_{11}^{-1}A_{12}$ and  $\tilde S = \tilde A_{22}-\tilde A_{12}^T \,  \tilde A_{11}^{-1} \, \tilde A_{12}$
denote  the Schur complements of $A$ and $\tilde A$,  respectively.
Then
\[
{\tilde A_{11}}^{-1}=(A_{11}+\tilde E_{11})^{-1}= A_{11}^{-1}-A_{11}^{-1} \tilde  E_{11} A_{11}^{-1}+{\cal O} (\normtwo{\tilde E_{11}}^2)
\]
and
\[
\tilde S=S +\tilde E_{22}+ W^T \tilde E_{11}W -(\tilde E_{12}^TW+W^T \tilde E_{12})+ {\cal O} (\normtwo{\tilde E}^2),
\]
where $W=A_{11}^{-1} A_{12}$.
\end{lemma}

\medskip
\begin{remark}\label{UwagaB1}
Under the  hypotheses  of Lemma \ref{bound_Schur} we have the bounds
\begin{equation}\label{bound_invA11}
\normtwo{{\tilde A_{11}}^{-1}- A_{11}^{-1}} \leq \normtwo{ A_{11}^{-1}}^2  \normtwo{\tilde E_{11}}+{\cal O} (\normtwo{\tilde E_{11}}^2)
\end{equation}
and
\begin{equation}\label{bound_tildeS}
\normtwo{\tilde S- S}  \leq   \normtwo{\tilde E_{22}}+ \normtwo{W}^2 \normtwo{\tilde E_{11}} + 2  \normtwo{W}  \normtwo{\tilde E_{12}} +  {\cal O} (\normtwo{\tilde E}^2).
\end{equation}
\end{remark}

\medskip

The next question is: {\bf  What can be said on the loss of symplecticity of  the matrices $L_1$ and $L_2$ from  Lemma \ref{lemacik}?  }
We recall that these matrices have the following forms
\begin{equation}\label{formsL1L2}
L_1=\left(
\begin{array}{cc}
 L_{11} &   0 \\
 L_{21} &  L_{11}^{-T}
\end{array}
\right),
\quad L_2=\left(
\begin{array}{cc}
 L_{11} &   0 \\
 L_{21} &  U_{22}
\end{array}
\right),
\end{equation}
where  $A_{11}=L_{11} L_{11}^T$ (Cholesky factorization), $L_{21}= (L_{11}^{-1} A_{12})^T$ and  $S=U_{22} U_{22}^T$ (Reverse Cholesky factorization of the Schur complement $S$).

\medskip
\begin{lemma}\label{sympfact1}
Let  $A \in \mathbb R^{2n \times 2n}$  be a symmetric positive definite  matrix   partitioned as in  (\ref{blockA}).
Let  $L_1, L_2  \in \mathbb R^{2n \times 2n}$   be given   in  (\ref{formsL1L2}).
Then we have
\begin{equation}\label{fact1}
\Omega(L_1) =\left(
\begin{array}{cc}
L_{11}^TL_{21}-L_{21}^TL_{11}  &  0 \\
0   &  0
\end{array}
\right),
\end{equation}
\begin{equation}
\label{fact1b}
\Omega(L_2) = \Omega(L_1) + \left(
\begin{array}{cc}
 0  & L_{11}^T U_{22}-I \\
-(L_{11}^T U_{22}-I)^T   & 0
\end{array}
\right).
\end{equation}
 \end{lemma}

\noindent
\begin{proof} 
It follows immediately from the formulae $\Omega(L_1)=L_1^TJL_1-J$ and $\Omega(L_2)=L_2^TJL_2-J$ for the matrices $L_1$ and $L_2$ from (\ref{formsL1L2}).
\end{proof}

\medskip
\begin{theorem}\label{twS1}
Under the  hypotheses  of Lemma \ref{sympfact1} we get the bounds
\begin{equation}\label{fact2}
\normtwo{\Omega(L_1)} \leq \normtwo{A_{11}^{-1}} \normtwo{\Omega(A)},
\end{equation}
and
\begin{equation}\label{fact4}
\normtwo{\Omega(L_1)}\leq  \normtwo{\Omega(L_2)} \leq  \normtwo{\Omega(L_1)}+ \normtwo{L_{11}^T U_{22}-I}.
 \end{equation}
Moreover, if
\begin{equation}\label{warunek1}
\rho(A_{11} (S-A_{11}^{-1}))\leq 1/2
\end{equation}
then
\begin{equation}\label{fact6}
\normtwo{L_{11}^T U_{22}-I} \leq {(2n)}^{1/2} \rho(A_{11} (S-A_{11}^{-1})).
\end{equation}
\end{theorem}

\noindent
\begin{proof} From (\ref{fact1})  and (\ref{Omega11})  it follows that  $L_{11}^TL_{21}-L_{21}^TL_{11}=L_{11}^{-1}\Omega_{11} L_{11}^{-T}$, hence  
$\normtwo{\Omega(L_1)} \leq \normtwo{L_{11}^{-1}}^2 \normtwo{\Omega_{11}}$.
From $\normtwo{\Omega_{11}} \leq  \normtwo{\Omega(A)}$ and the identity $\normtwo{L_{11}^{-1}}^2 = \normtwo{A_{11}^{-1}}$ we obtain  (\ref{fact2}).
Clearly, the result (\ref{fact4}) follows immediately  from the property of 2-norm.
From $U_{22} U_{22}^T =S = A_{11}^{-1} + (S-A_{11}^{-1}) $ we have the identity
\[
(L_{11}^T U_{22}) (L_{11}^T U_{22})^T = I + L_{11}^T (S-A_{11}^{-1}) L_{11}.
\]
Note that the matrix $L_{11}^T U_{22}$ is upper triangular with positive diagonal entries. Thus, we get the unique  Reverse factorization  of the matrix $I + L_{11}^T (S-A_{11}^{-1}) L_{11}$. 
To prove (\ref{fact6}), we apply Theorem \ref{ReChol1}. Here the matrix $I$ takes on the role of the matrix $A$,  $U=I$   and $E= L_{11}^T (S-A_{11}^{-1}) L_{11}$. 
Notice that $E$ is symmetric, hence  $\normtwo{E}=\rho(E)$. Since $A_{11}=L_{11} L_{11}^T$, we get  $E=L_{11}^{-1} (A_{11} (S-A_{11}^{-1})) L_{11}$, hence 
$\normtwo{E}=\rho(E)=\rho(A_{11} (S-A_{11}^{-1}))$. 
Under the assumption \ref{warunek1},  Theorem \ref{ReChol1}   says for $p=2$ that
\[
||L_{11}^{T} U_{22}  - I||_{F} \leq  2^{-1/2}  \frac{||E||_{F}}{1-\normtwo{E}}.
\]

This together with the inequalities $\normtwo{L_{11}^{T} U_{22}  - I} \leq ||L_{11}^{T} U_{22}  - I||_F$ and   $||E||_{F} \leq \sqrt{n} \normtwo{E}$ leads to (\ref{fact6}).
The proof is complete.
\end{proof}

\medskip
\begin{remark}\label{AAA}
{\bf  Theorem \ref{twS1}  says that the departure of  the matrix $L_1$ from symplecticity is not greater than for  the matrix $L_2$.

If  we replace (\ref{warunek1}) by 
\begin{equation}\label{warunek2}
\normtwo{A_{11}} \normtwo{S-A_{11}^{-1}}\leq 1/2
\end{equation}
then we have a bound
\begin{equation}\label{fact7}
\normtwo{L_{11}^T U_{22}-I} \leq {(2n)}^{1/2} \normtwo{A_{11}} \normtwo{S-A_{11}^{-1}}.
\end{equation}

Since $\rho(X)\leq \normtwo{X}$ for any matrix $X(2n \times 2n)$, we see that the result stated in (\ref{fact6}) is stronger than given in    (\ref{fact7}).
}
\end{remark}

\medskip
If  $A  \in \mathbb R^{2n \times 2n}$ is symplectic then  its inverse equals $A^{-1}=J^TA^TJ$, and it is also symplectic.
Then the condition number of $A$ equals $\kappa_2(A)=\normtwo{A}^2$.
 \begin{lemma} [{\cite[Lemma 8]{Maks}}]
\label{Maks1}
Let $A \in \mathbb R^{2n \times 2n}$ satisfy  $\normtwo{\Omega(A)}<1$.
Then $A$ is nonsingular and we have
\[
\kappa_2 (A) \leq \frac{\normtwo{A}^2}{1-\normtwo{\Omega(A)}}.
\]
\end{lemma}

\noindent

In the following section, we will prove  that Algorithm $W_2$  is numerically stable, like a standard Cholesky factorization.
This algorithm  can be applied to a broader class of symmetric positive definite matrices $A$  that are not necessarily  symplectic.
However, if $A$ is additionally symplectic, then we get the symplectic  factorization $A=LL^T$, with symplectic factor  $L$.
We will also show that Algorithm $W_1$ can be less stable than Algorithm $W_2$. In particular, if we apply Algorithm $W_1$
to a matrix $A$ that is only symmetric positive definite, the factorization error depends on the quantity $\normtwo{A_{11}^{-1}-S}$,
measuring the distance of $A$ from symplecticity, and on the condition number $\kappa(A_{11})=\normtwo{A_{11}^{-1}} \cdot \normtwo{A_{11}}$ of the principal submatrix $A_{11}$.

\medskip

\section{Rounding error analysis}\label{Error analysis}

Let $u$ denote   machine precision. To simplify expressions  in rounding error analysis we  use  $\gamma_n(\eps)= \frac{n \eps}{1 - n \eps}$ for $n \eps<1$, where
$\eps$ denotes  machine precision.
First we recall some basics lemmas and theorems.

\begin{lemma} \label{Lemat0}
Let $C=A B$, where $A, B \in \mathbb R^{n \times n}$. Then the computed product $\tilde C$ of $A$ and $B$ by standard matrix multiplication algorithm satisfies
\[
\tilde C =  C+ \Delta C, \hskip3ex  |\Delta C| \leq \gamma_n(\eps) |A| \, |B|.
\]
\end{lemma}
We will also use some rounding error results on the solution of triangular systems, see \cite[pp.\ 142, 262--263]{higham:2002}.
\begin{lemma} \label{Lemat1}
Let $L \in \mathbb R^{n \times n}$ be a nonsingular lower triangular matrix and let $Y \in \mathbb R^{n \times n}$.
Then the  solution $\tilde X \in \mathbb R^{n \times n}$ to $LX =Y$  computed  by the forward  substitution satisfies
\[
L \tilde X = Y+ \Delta Y, \hskip3ex |\Delta Y| \leq \gamma_n(\eps) |L| \, |\tilde X|,
\hskip1ex \text{and} \hskip2ex \normtwo{\Delta Y} \leq n \gamma_n(\eps)  \normtwo{L} \normtwo{\tilde X}.
\]
\end{lemma}
In the following we also recall some results on several factorizations of the symmetric positive definite  matrix $A$.
\begin{lemma}\label{Lemat2}
Let $A \in \mathbb R^{n \times n}$ be symmetric positive definite,  $\Delta A  \in \mathbb R^{n \times n}$ be symmetrix
and  $B \in \mathbb R^{n \times n}$.
Let $A+\Delta A=BB^T$, where $|\Delta A| \leq  \delta |B| \, |B^T|$  for some constant $\delta>0$ such that $n \delta<1$.
Then $\normtwo{\Delta A} \leq \gamma_n(\delta)  \normtwo{A}$,
where  $\gamma_n(\delta)= \frac{n \delta}{1 - n \delta}$.
\end{lemma}
\noindent
\begin{proof}
We use the same reasoning as in \cite{higham:2002}, pp. 198.  The key inequality is
\[
\normtwo{|B||B^T|} = \normtwo{|B|}^2 \leq n \normtwo{B}^2= n \normtwo{BB^T}=n \normtwo{A+\Delta A}.
\]
This together with the inequality $\normtwo{\Delta A} \leq \delta \normtwo{|B||B^T|}$ leads to $\normtwo{\Delta A} \leq n \delta  \normtwo{A+\Delta A}$.
Note that $\normtwo{A+\Delta A} \leq \normtwo{A}+\normtwo{\Delta A}$, hence the lemma follows.
\end{proof}

\medskip
\noindent
The next result is on numerical properties of  the Cholesky factorization, see  Theorem 10.3  \cite[pp. 197]{higham:2002}.
We apply  Lemma \ref{Lemat2} to estimate  $\normtwo{\Delta A}$ and to establish Theorems \ref{Cholesky} and   \ref{ReverseCholesky}.

\begin{theorem} \label{Cholesky}
If Cholesky factorization applied to the symmetric and positive definite matrix  $A \in \mathbb R^{n \times n}$
runs to completion, then  the computed  lower triangular factor $\tilde L$   satisfies
\[
A+\Delta A= \tilde L {\tilde L}^T, \hskip3ex |\Delta A| \leq \gamma_{n+1}(\eps) |\tilde L| \, |{\tilde L}^T|,
\]
and
\[ \normtwo{\Delta A} \leq n \gamma_{n+1}(\eps)(1-n \gamma_{n+1}(\eps))^{-1} \,  \normtwo{A} \leq  2n \gamma_{n+1}(\eps) \, \normtwo{A},
\]
provided that $2 n \gamma_{n+1}(\eps) <1$.
\end{theorem}
Clearly, the following analogue of  Theorem  \ref{Cholesky} holds.
\begin{theorem} \label{ReverseCholesky}
If Reverse Cholesky factorization applied to the symmetric and positive definite matrix  $A \in \mathbb R^{n \times n}$
runs to completion then  the computed  upper triangular factor $\tilde U$   satisfies
\[
A+\Delta A= \tilde U {\tilde U}^T, \hskip3ex |\Delta A| \leq \gamma_{n+1}(\eps)  |\tilde U| \, |{\tilde U}^T|,
\]
and
\[ \normtwo{\Delta A} \leq n \gamma_{n+1}(\eps)(1-n \gamma_{n+1}(\eps))^{-1} \,  \normtwo{A} \leq  2n \gamma_{n+1}(\eps) \, \normtwo{A},
\]
provided that $2 n \gamma_{n+1}(\eps)<1$.
\end{theorem}

\noindent
Based on the previous results, the following  lemma can be given for a general $LL^T$ factorization
of the symmetric positive definite matrix  partitioned  as in (\ref{blockA}).
\begin{lemma}\label{Lemat3}
Let $A \in \mathbb R^{2n \times 2n}$ be a symmetric positive definite  matrix partitioned  as in (\ref{blockA}).
Let $\tilde L  \in \mathbb R^{2n \times 2n}$  be a nonsingular  block lower triangular matrix.
Define $E=\tilde L {\tilde L}^T - A$.  Partition $\tilde L$ and $E$ conformally with $J$ as follows
\[
\tilde L=\left(
\begin{array}{cc}
 \tilde L_{11} &   0 \\
\tilde  L_{21} & \tilde  L_{22}
\end{array}
\right),
\quad
E=\left(
\begin{array}{cc}
E_{11} &   E_{12} \\
E_{12}^T & E_{22}
\end{array}
\right).
\]

Then  the blocks of the matrix $E$  satisfy the following identities
\begin{eqnarray}\label{E11}
E_{11} & = & \tilde L_{11} \tilde L_{11}^T-A_{11},
\\
\label{E12}
E_{12} & = & \tilde L_{11} \tilde L_{21}^T-A_{12},
\\
\label{E22}
E_{22} & = &  \tilde L_{22} \tilde L_{22}^T-(A_{22}-\tilde L_{21} \tilde L_{21}^T).
\end{eqnarray}
\end{lemma}

\begin{corollary}\label{corollaryB1}
Considering  Lemma  \ref{Lemat3}, if  $\tilde L_{11}$ and $\tilde L_{21}$ are  computed by Algorithm $W_1$ or $W_2$, i.e.
by the Cholesky factorization of $A_{11}$ and by consequent forward substitution applied to $\tilde L_{11}$,
then from Theorem \ref{Cholesky} it follows that
\begin{equation}\label{boundE11}
\tilde A_{11}= A_{11}+E_{11}=\tilde L_{11} \tilde L_{11}^T, \quad  |E_{11}|  \leq \gamma_{n+1}(\eps) |\tilde L_{11}|\,|\tilde L_{11}^T|.
\end{equation}
Moreover, from Lemma \ref{Lemat1} we have
\begin{equation}\label{boundE12}
\tilde A_{12}= A_{12}+E_{12}=\tilde L_{11} \tilde L_{21}^T, \quad |E_{12}|  \leq \gamma_{n}(\eps)  |\tilde L_{11}| \, |\tilde L_{21}^T|.
\end{equation}
\end{corollary}

\medskip
\noindent
In the following we obtain results for Algorithm $W_2$
that are similar  to those for Cholesky factorization, see Theorem \ref{Cholesky}.
\begin{theorem}\label{ThmW2}
If Algorithm $W_2$ applied to the symmetric and positive definite matrix  $A \in \mathbb R^{2n \times 2n}$
runs to completion,  then  the computed  block triangular factor $\tilde L_2  \in \mathbb R^{2n \times 2n}$   satisfies
\begin{equation}\label{boundW2}
A+E= \tilde L_2 {\tilde L_2}^T, \hskip3ex  |E| \leq \eps |A|+ \gamma_{n+1}(\eps)  |\tilde L_2|  |\tilde L_2^T|.
\end{equation}
Moreover,  the factorization error satisfies the normwise bound
\begin{equation}\label{NormdecW2}
\normtwo{A-\tilde L_2 {\tilde L_2}^T}=\normtwo{E} \leq  4n \gamma_{n+2}(\eps)  \normtwo{A},
\end{equation}
provided that $4n \gamma_{n+2}(u)<1$.
\end{theorem}
\noindent
\begin{proof}
Let $\tilde S$ be the computed Schur complement $S=A_{22}-L_{21}  L_{21}^T$ via  standard method for product of two matrices.
By Lemma \ref{Lemat0},  we have
\begin{equation}\label{tildeS}
\tilde S+\Delta \tilde S= A_{22}-\tilde L_{21} \, {\tilde L_{21}}^T, \quad |\Delta \tilde S| \leq \eps |A_{22}|  + \gamma_{n+1}(\eps) |\tilde L_{21}| \, |\tilde L_{21}^T|.
\end{equation}
From Theorem \ref{ReverseCholesky} it follows that
\begin{equation}\label{tilde L22}
\tilde S +\Delta  S= \tilde L_{22} \, {\tilde L_{22}}^T, \quad |\Delta S| \leq \gamma_{n+1}(\eps) |\tilde L_{22}| \, |\tilde L_{22}^T|.
\end{equation}
This together with (\ref{E22}) leads to
\[
E_{22}=(\tilde S +\Delta  S)-(\tilde S +\Delta \tilde S)= \Delta  S - \Delta \tilde S.
\]
Then  $|E_{22}| \leq  |\Delta \tilde S| + |\Delta S|$,   so making use of  (\ref{tildeS})-(\ref{tilde L22}) we conclude that
\[
|E_{22}| \leq \eps |A_{22}|  + \gamma_{n+1}( \eps) |\tilde L_{21}| |\tilde L_{21}^T| +  \gamma_{n+1}(\eps) |\tilde L_{22}| \, |\tilde L_{22}^T|.
\]
Consequently,  this together with (\ref{boundE11})-(\ref{boundE12})  leads to
\[
|E| \leq \eps
\left(
\begin{array}{cc}
0 &  0  \\
0 & |A_{22}|
\end{array}
\right)
+ \gamma_{n+1}(\eps)
\left(
\begin{array}{cc}
|\tilde L_{11}| |\tilde L_{11}^T| &  |\tilde L_{11}|  |\tilde L_{21}^T|  \\
|\tilde  L_{21}|  |\tilde L_{11}^T|  & |\tilde L_{21}| |\tilde L_{21}^T| + |\tilde L_{22}| |\tilde L_{22}^T|
\end{array}
\right).
\]
which we can rewrite as
\[
|E| \leq \eps
\left(
\begin{array}{cc}
0 &  0  \\
0 & |A_{22}|
\end{array}
\right)
+\gamma_{n+1}(\eps) |\tilde L_2| \,  |\tilde L_2^T| \leq \eps |A| + \gamma_{n+1}(\eps) |\tilde L_2| \,  |\tilde L_2^T|.
\]
This establishes the formula (\ref{boundW2}).
To prove (\ref{NormdecW2}), we apply the same idea as in  Lemma \ref{Lemat2}.
Taking  norms  in  (\ref{boundW2}),    we  obtain
\[
\normtwo{E} \leq \eps \normtwo{|A|}+ \gamma_{n+1}(\eps)  \normtwo{|\tilde L_2| \,  |\tilde L_2^T|}
\leq \sqrt{2n} \eps \normtwo{A}+  2n  \gamma_{n+1}(\eps)  \normtwo{\tilde L_2}^2.
\]
Since $\tilde L_2 \,  \tilde L_2^T=A+E$, hence  $\normtwo{\tilde L_2}^2=\normtwo{\tilde L_2 \,  \tilde L_2^T}\leq \normtwo{A}+ \normtwo{E}$.
Thus,
\[
\normtwo{E} \leq  (\sqrt{2n} \eps +2n  \gamma_{n+1}(\eps))  \normtwo{A}+ 2n  \gamma_{n+1}(\eps) \normtwo{E}.
\]
Now we use  the bounds  $\sqrt{2n}\leq 2n$ and $u+\gamma_{n+1}(\eps) \leq \gamma_{n+2}(\eps)$. Then
\[
\normtwo{E} \leq  2n  \gamma_{n+2}(\eps)  \normtwo{A}+ 2n  \gamma_{n+1}(\eps) \normtwo{E},
\]
hence
$\normtwo{E} \leq  2n  \gamma_{n+2}(\eps) (1- 2n  \gamma_{n+1}(\eps))^{-1} \,  \normtwo{A} \leq 4n  \gamma_{n+2}(\eps)  \,  \normtwo{A}$.
This gives the formula (\ref{NormdecW2}), and the proof is complete.
 \end{proof}

\medskip
\noindent
Theorem  \ref{ThmW2}  states  that  Algorithm $W_2$ is numerically stable for a  class of symmetric positive definite matrices $A \in \mathbb R^{2n \times 2n}$,  producing the computed factor $\tilde L_2$  such that $\normtwo{A-\tilde L_2 {\tilde L_2}^T}=  {\cal O}(\eps \normtwo{A})$.

\medskip
The next question is:  {\em What can be said on numerical stability of  Algorithm $W_1$?}
Let the matrix $\tilde L_1=\left(
\begin{array}{cc}
\tilde L_{11} &  0  \\
\tilde L_{21} & \tilde L_{22}
\end{array}
\right)$ be   computed from Algorithm $W_1$. We recall that  the matrices  $\tilde L_{11}$ and $\tilde L_{21}$ are exactly the same as for Algorithm $W_2$.
However, in contrast to  Algorithm $W_2$, here   $\tilde L_{22} \approx \tilde L_{11}^{-T}$ is obtained by inverting lower triangular matrix $\tilde L_{11}$.
The following lemma gives the bounds for the computed matrix $\tilde L_{22}$ by Algorithm $W_1$.
\medskip
\begin{lemma}\label{L22_W1}
Let $A \in \mathbb R^{2n \times 2n}$ be symmetric and positive definite. Assume that Algorithm $W_1$ applied to $A$ runs to completion,
producing  the block lower triangular matrix  $\tilde L_1 \in \mathbb R^{2n \times 2n}$. Let  $\tilde A_{11}=A_{11}+E_{11}$ be  defined in  (\ref{boundE11}).
Then  $\tilde L_{22}$ satisfies
\begin{equation}\label{tildeL22}
\tilde L_{22} \tilde L_{22}^T=\tilde  A_{11}^{-1}+ G_{22},
\end{equation}
where
\begin{equation}\label{G22}
\normtwo{G_{22}} \leq 2 n \gamma_n(\eps)  \normtwo{\tilde A_{11}^{-1}} \, {\kappa_2(\tilde A_{11})^{1/2}}+ {\cal O} (\eps^2).
\end{equation}
\end{lemma}
\noindent
\begin{proof}
From  Lemma \ref{Lemat1} it follows that
\begin{equation}\label{invL11}
\tilde L_{11} \tilde L_{22}^T = I+ F_{22}, \hskip3ex  \normtwo{F_{22}}  \leq n \gamma_n(\eps) \normtwo{\tilde L_{11}} \, \normtwo{\tilde L_{22}^T}.
\end{equation}
Thus,
$ \tilde L_{22} \tilde L_{22}^T=(I+F_{22}^T) (\tilde L_{11}^{-T}\, \tilde L_{11}^{-1}) (I+F_{22})$.
From   (\ref{boundE11}) we have   $\tilde  A_{11}^{-1}=(\tilde L_{11} \tilde L_{11}^T)^{-1}= \tilde L_{11}^{-T}\, \tilde L_{11}^{-1}$, hence
the above equation can be written as follows
\[
\tilde L_{22} \tilde L_{22}^T=(I+F_{22}^T) \tilde A_{11}^{-1} (I+F_{22}) =\tilde  A_{11}^{-1}+ G_{22},
\]
where $G_{22}=\tilde  A_{11}^{-1} F_{22} +(\tilde  A_{11}^{-1} F_{22})^T+ F_{22}^T \tilde A_{11}^{-1} F_{22}$.
Taking norms, we obtain  the following bound
\[
\normtwo{G_{22}} \leq 2 \normtwo{\tilde A_{11}^{-1} F_{22}} +  \normtwo{\tilde A_{11}^{-1}} \normtwo{F_{22}}^2
\leq 2 \normtwo{\tilde A_{11}^{-1}} \normtwo{F_{22}} + {\cal O} (\eps^2).
\]
This together with (\ref{invL11}) implies  the bound
\[
\normtwo{G_{22}} \leq 2 n \gamma_n(\eps) \normtwo{\tilde A_{11}^{-1}} \normtwo{\tilde L_{11}} \, \normtwo{\tilde L_{22}^T} + {\cal O} (\eps^2).
\]
From  (\ref{invL11}) we get    $\tilde L_{22}^T = \tilde L_{11}^{-1} + \tilde L_{11}^{-1} F_{22}$. This together with
$\kappa_2(\tilde A_{11})= \normtwo{\tilde L_{11}}^2 \, \normtwo{\tilde L_{11}^{-1}}^2$ leads to
$
\normtwo{G_{22}} \leq 2 n \gamma_n(\eps)  \normtwo{\tilde A_{11}^{-1}} \, {\kappa_2(\tilde A_{11})^{1/2}}+ {\cal O} (\eps^2)$.
This completes the proof.
\end{proof}

\medskip

\begin{lemma}\label{LematW1}
Under the  hypotheses of   Lemma  \ref{L22_W1}, let $\tilde A_{11}=A_{11}+E_{11}$  and  $\tilde A_{12}=A_{12}+E_{12}$
be given as in (\ref{boundE11})-(\ref{boundE12}).
Assume that  the matrix
\begin{equation}\label{EhatA}
\tilde A=
\left(
\begin{array}{cc}
\tilde A_{11}   &  \tilde A_{12}  \\
{\tilde A_{12}}^T  & A_{22}
\end{array}
\right)
\end{equation}
is  positive definite.  Then the factorization error satisfies
\begin{equation}\label{NormdecW1}
\normtwo{\tilde A-\tilde L_1 {\tilde L_1}^T} \leq \normtwo{\tilde A_{11}^{-1}  - \tilde S} +  2n \gamma_{n}(\eps)  \normtwo{\tilde A_{11}^{-1}} \, {\kappa_2(\tilde A_{11})^{1/2}}+ {\cal O} (\eps^2),
\end{equation}
where $\tilde S= A_{22}-\tilde A_{12}^T \,  \tilde A_{11}^{-1} \, \tilde A_{12}$ is the Schur complement of $\tilde A_{11}$  for the matrix $\tilde A$.
\end{lemma}
\noindent
\begin{proof}
Notice that  we have $\tilde A=A+ \tilde E$, where
\begin{equation}\label{tildeE}
\tilde E
=\left(
\begin{array}{cc}
E_{11}   &  E_{12}  \\
E_{12}^T & 0
\end{array}
\right).
\end{equation}
From Lemma  \ref{Lemat3} it follows that  the matrix $\tilde L_1=\left(
\begin{array}{cc}
\tilde L_{11} &  0  \\
\tilde L_{21} & \tilde L_{22}
\end{array}
\right)$
satisfies
\begin{equation}\label{hatAE}
\tilde L_1 {\tilde L_1}^T- \tilde A
=\left(
\begin{array}{cc}
0   &  0  \\
0 & E_{22}
\end{array}
\right),
\end{equation}
where $E_{22}$ is defined in (\ref{E22})  as
\begin{equation}\label{W1E22}
E_{22}=\tilde L_{22} \tilde L_{22}^T- (A_{22}-\tilde L_{21} \tilde L_{21}^T).
\end{equation}
From (\ref{boundE11})-(\ref{boundE12}) we get the identities
$\tilde L_{21}^T= \tilde L_{11} ^{-1} \tilde A_{12}$ and   $\tilde A_{11} ^{-1}=\tilde L_{11} ^{-T}   \tilde L_{11} ^{-1}$.
We see  that
$A_{22}- \tilde L_{21} \tilde L_{21}^T$
is in fact  the Schur complement $\tilde S= A_{22}-\tilde A_{12}^T \,  \tilde A_{11}^{-1} \, \tilde A_{12}$  of $\tilde A_{11}$ for  the matrix $\tilde A$.
From  Lemma  \ref{L22_W1} and  (\ref{tildeE})-(\ref{W1E22}) it follows that
\[
\normtwo{\tilde  A- \tilde L_1 {\tilde L_1}^T}= \normtwo{E_{22}}= \normtwo{(\tilde  A_{11}^{-1}+G_{22})- (A_{22} - \tilde A_{12}^T \,  \tilde A_{11}^{-1} \, \tilde A_{12})} \]
\[
=\normtwo{(\tilde A_{11}^{-1}- \tilde S) +G_{22}}\leq \normtwo{\tilde A_{11}^{-1}- \tilde S }+ \normtwo{G_{22}}.
\]
This together with (\ref{G22})  leads to (\ref{NormdecW1}). The proof is complete.
\end{proof}

\medskip
\begin{corollary}\label{Corollary A1}
Under the  hypotheses of   Lemma  \ref{LematW1}, we have the bound
\begin{equation}\label{closeA}
\normtwo{\tilde A-A} = \normtwo{\tilde E} \leq 2 n \gamma_{n+1}(\eps) \normtwo{A} + {\cal O} (\eps^2),
\end{equation}
provided that $2 n \gamma_{n+1}(\eps)<1$.

\noindent
\begin{proof}
From (\ref{tildeE}) and by the property of 2-norm we get
\[
\normtwo{\tilde E} \leq \normtwo{E_{11}}+\normtwo{E_{12}}.
\]
From (\ref{boundE11})-(\ref{boundE12}) it follows that
\[
\normtwo{E_{11}} \leq n \gamma_{n+1}(u) \normtwo{\tilde L_{11}}^2 \leq n \gamma_{n+1}(u) \normtwo{L_{11}}^2+{\cal O} (\eps^2)
\]
and
\[
\normtwo{E_{12}} \leq n \gamma_{n}(u) \normtwo{\tilde L_{11}} \, \normtwo{\tilde L_{21}} \leq n \gamma_{n}(u) \normtwo{L_{11}} \, \normtwo{L_{21}}+{\cal O} (\eps^2).
\]
Notice that $\normtwo{L_{21}} \leq \normtwo{L_2}=\normtwo{A}^{1/2}$, where $L_2$ is given  in Theorem \ref{AlgW2}.
Since $\normtwo{A_{11}} \leq \normtwo{A}$, we obtain the inequalities
\begin{equation}\label{boundForE11}
\normtwo{E_{11}} \leq n \gamma_{n+1}(u)   \normtwo{A_{11}}+{\cal O} (\eps^2) \leq n \gamma_{n+1}(u)   \normtwo{A}+{\cal O} (\eps^2)
\end{equation}
and
\begin{equation}\label{boundForE12}
\normtwo{E_{12}} \leq n \gamma_{n}(\eps) \normtwo{A_{11}}^{1/2}  \normtwo{A}^{1/2}+{\cal O} (\eps^2).
\end{equation}
This leads to
\[
\normtwo{\tilde A-A} = \normtwo{\tilde E} \leq \normtwo{E_{11}}+\normtwo{E_{12}} \leq 2 n \gamma_{n+1}(\eps) \normtwo{A} + {\cal O} (\eps^2),
\]
which proves (\ref{closeA}).
\end{proof}
\end{corollary}

\medskip

The bound  (\ref{closeA}) confirms  that  $\tilde A$  is close to the matrix $A$.
It is  well-known that if $\normtwo{A^{-1}} \, \normtwo{\tilde E}<1$ then the matrix $\tilde A$ is also positive definite.

\medskip

The task is now to estimate  the quantity $\normtwo{\tilde A_{11}^{-1}  - \tilde S}$  by $\normtwo{A_{11}^{-1}-S}$, where $S$ is the  Schur complement of $A_{11}$ for the matrix $A$.
We resolve this problem with help of Lemma \ref{bound_Schur} for the matrix $\tilde E$ defined in (\ref{tildeE}).

\medskip
\begin{theorem}\label{NewThmW1}
Under the  hypotheses of Lemma \ref{LematW1},  we have the following bounds
\begin{equation}\label{bound invA11b}
\normtwo{{\tilde A_{11}}^{-1}- A_{11}^{-1}} \leq n \gamma_{n+1}(\eps) \kappa_2(A_{11}) \normtwo{ A_{11}^{-1}}+{\cal O} (\eps^2),
\end{equation}
\begin{equation}\label{new_boundS}
\normtwo{\tilde S- S}\leq   3n\gamma_{n+1}(\eps) \kappa_2(A_{11}) \,  \normtwo{A}+{\cal O} (\eps^2),
\end{equation}
and
\begin{equation}\label{bound_tA}
\normtwo{\tilde A_{11}^{-1}- \tilde S} \leq \normtwo{A_{11}^{-1}-S}+n \gamma_{n+1}(\eps) \kappa_2(A_{11}) (\normtwo{ A_{11}^{-1}}+ 3\normtwo{A}) +{\cal O} (\eps^2).
\end{equation}
Moreover, the factorization error in  Algorithm $W_1$ satisfies the bound
\begin{eqnarray}
\normtwo{A-\tilde L_1 {\tilde L_1}^T} & \leq & \normtwo{A_{11}^{-1}-S} +  2n \gamma_{n+1}(\eps)  \normtwo{A}
\nonumber  \\
& + & 3 n\gamma_{n+1}(\eps) \kappa_2(A_{11}) (\normtwo{A}+ \normtwo{A_{11}^{-1}})  + {\cal O} (\eps^2).\label{SecondW1}
\end{eqnarray}
\end{theorem}
\noindent
\begin{proof}
From (\ref{bound_invA11}) it follows that
\[
\normtwo{\tilde A_{11}^{-1}- A_{11}^{-1}}\leq \normtwo{E_{11}} \normtwo{ A_{11}^{-1}}^2 +{\cal O} (\eps^2).
\]
Since  (\ref{boundForE11}) says that $\normtwo{E_{11}} \leq n \gamma_{n+1}(u)   \normtwo{A_{11}}+{\cal O} (\eps^2)$
and $\kappa_2(A_{11})=\normtwo{A_{11}}\, \normtwo{A_{11}^{-1}}$, we  obtain   (\ref{bound invA11b}).
Next, from (\ref{bound_tildeS}) we get
\[
\normtwo{\tilde S -S}\leq \normtwo{E_{11}} \normtwo{W}^2 + 2 \normtwo{E_{12}} \normtwo{W}+{\cal O} (\eps^2).
\]
This together with (\ref{boundForE11})-(\ref{boundForE12}) leads to
\[
\normtwo{\tilde S- S}\leq   n\gamma_{n+1}(\eps) (\normtwo{A_{11}}\, \normtwo{W}^2 + 2 \normtwo{A_{11}}^{1/2} \normtwo{A}^{1/2} \normtwo{W})+{\cal O} (\eps^2),
\]
where $W=A_{11}^{-1} A_{12}$.
From Lemma \ref{boundW} and the inequality $\normtwo{A_{22}} \leq \normtwo{A}$ it follows that $\normtwo{W}^2 \leq \normtwo{ A_{11}^{-1}}\, \normtwo{A}$, hence
\[
\normtwo{\tilde S- S}\leq   n \gamma_{n+1}(\eps) (\kappa_2(A_{11}) + 2 \kappa_2(A_{11})^{1/2})  \, \normtwo{A} +{\cal O} (\eps^2),
\]
which implies   (\ref{new_boundS}).
Taking norms in the expression
\[
\tilde A_{11}^{-1}- \tilde S= (A_{11}^{-1}-S) + (\tilde A_{11}^{-1}-A_{11}^{-1}) - (\tilde S- S)
\]
leads to
$\normtwo{\tilde A_{11}^{-1}- \tilde S} \leq \normtwo{A_{11}^{-1}- S}+\normtwo{\tilde A_{11}^{-1}- A_{11}^{-1}} +  \normtwo{\tilde S- S}$.
Combining this with (\ref{bound invA11b})-(\ref{new_boundS}) yields (\ref{bound_tA}).
Similarly, we see  that  the formula
\begin{equation}
\label{FirstW1}
\normtwo{\tilde A-\tilde L_1 {\tilde L_1}^T}  \leq  \normtwo{A_{11}^{-1}-S} + 3n \gamma_{n+1}(\eps) \kappa_2(A_{11}) (\normtwo{A_{11}^{-1}}+\normtwo{A})+ {\cal O} (\eps^2)
\end{equation}
follows  immediately from (\ref{NormdecW1}) and (\ref{bound_tA}).
Writing  $A-\tilde L_1 {\tilde L_1}^T=(\tilde  A-\tilde L_1 {\tilde L_1}^T)- (\tilde A-A)$, and taking the norms, we have
\[
\normtwo{A-\tilde L_1 {\tilde L_1}^T} \leq \normtwo{\tilde  A-\tilde L_1 {\tilde L_1}^T}+ \normtwo{\tilde A-A}.
\]
Since $\tilde A$ is close to $A$, see (\ref{closeA}) in Corollary \ref{Corollary A1}, we conclude that
\[
\normtwo{A-\tilde L_1 {\tilde L_1}^T} \leq  \normtwo{\tilde  A-\tilde L_1 {\tilde L_1}^T}+ 2n \gamma_{n+1}(\eps) \normtwo{A}+{\cal O} (\eps^2).
\]
This together with (\ref{FirstW1}) establishes the formula (\ref{SecondW1}).
The proof is complete.
\end{proof}

\medskip

\begin{remark}\label{RemarkB}
Let us mention some  consequences of Theorem  \ref{NewThmW1}.  How to relate $\normtwo{A_{11}^{-1}}$ to $\normtwo{A}$?
From  (\ref{sympL})  and the property of 2-norm it follows that $\normtwo{L_{22}}=\normtwo{L_{11}^{-T}}\leq \normtwo{L_1}$, hence
$\normtwo{A_{11}^{-1}}=\normtwo{L_{22}}^2 \leq \normtwo{L_1}^2$.
Then  (\ref{Miro1}) implies
$\normtwo{L_1}^2 \leq   \normtwo{A} + \normtwo{\Delta A}= \normtwo{A} + \normtwo{A_{11}^{-1}-S}$.
Thus,
\begin{equation}\label{AA1}
\normtwo{A_{11}^{-1}}\leq \normtwo{A} + \normtwo{A_{11}^{-1}-S}.
\end{equation}
We conclude that (\ref{AA1}) together with (\ref{SecondW1}) leads to the factorization error bound for Algorithm $W_1$
\begin{eqnarray}
 \label{AA2}
\normtwo{A-\tilde L_1 {\tilde L_1}^T} & \leq  &  \normtwo{A_{11}^{-1}-S}  \left ( 1 + 3n \gamma_{n+1}(\eps) \kappa_2(A_{11})  \right )
\nonumber  \\
 &   +  & 8n \gamma_{n+1}(\eps) \kappa_2(A_{11})  \normtwo{A}+{\cal O} (\eps^2).
\end{eqnarray}
In addition,
if  $\normtwo{A_{11}^{-1}-S}/\normtwo{A}={\cal O}(\eps)$ and the matrix $A_{11}$ is well-conditioned, i.e.,  $\kappa_2(A_{11})={\cal O} (1)$, then Algorithm $W_1$ is
numerically stable.
\end{remark}

\medskip

\section{Numerical tests}\label{tests}

This section  reports on  the results of numerical experiments that examine
the behavior  of considered algorithms.
All tests are performed in \textsl{MATLAB} ver. R2023b,  with the machine precision $\eps \approx 2.2 \cdot 10^{-16}$.
We consider several examples of symmetric positive definite matrices that are also symplectic or near-symplectic, where the  distance
to symplecticity of the matrix $A$ is measured by the  quantities
$\normtwo{A^TJA-J}$ and $\normtwo{A_{11}^{-1} - S}/\normtwo{A}$. We look at the numerical behavior
of Algorithm  $W_1$ and Algorithm $W_2$,  and report the
factorization error $\normtwo{A-L_1  L_1^T}/\normtwo{A}$ for Algorithm $W_1$ and the factorization error $\normtwo{A-L_2  L_2^T}/\normtwo{A}$
for Algorithm $W_2$, respectively.
{\bf We examine the bounds obtained  for the factorization error for both algorithms in terms of the quantity $\normtwo{A_{11}^{-1} - S}$, computed as described in  (\ref{computeL1})  in  Remark \ref{UwagaA2}. Moreover,
we display the condition number   $\kappa_2 (A_{11})$ and  $\normtwo{A}$. }In addition, we report the quantities  on the loss of symplecticity  of the matrices
$A^TJA-J$,  $L_1^TJL_1-J$  for  Algorithm $W_1$,
and $L_2^TJL_2-J$ for  Algorithm $W_2$.

\medskip

\begin{example}\label{example1}
First,
we  test Algorithms $W_1$ and $W_2$ on the matrix $A= S^T\, S$, where $S$ is considered in  \cite{Tam} (see also  \cite{Benzi}-\cite{Maks}):
\begin{equation}\label{S}
S=S(\theta)=\left(
\begin{array}{cccc}
\cosh{\theta}   &   \sinh{\theta} & 0 &   \sinh{\theta}\\
 \sinh{\theta} & \cosh{\theta} &  \sinh{\theta} & 0 \\
0 & 0 & \cosh{\theta}   &  - \sinh{\theta} \\
0 & 0 & -\sinh{\theta}   &   \cosh{\theta} \\
 \end{array}
\right),
\quad \theta \in  \mathbb R.
\end{equation}

{\bf Note that $S$ and $A$ are  symplectic in the theory  but their representation in floating point arithmetic make them to be only   near symplectic.
In fact, here $A={\tilde S}^T {\tilde S}$, where $\tilde S=fl(S)$. } 

\medskip
\begin{table}\label{tabelka1}
\caption{Results  for Example \ref{example1} with $A=S^TS$, where $S$ is defined in (\ref{S}).}
\begin{center}
\begin{tabular}{lccccc}
	\hline
	$\theta$ &  $3$ & $4$ & $6$ & $7$\\
	\hline
        $\kappa_2 (A)$ &   2.5380e+05  & 1.3881e+07 &  4.1389e+10 &  2.2601e+12  \\
        $\normtwo{A}$ &   5.0379e+02 &   3.7257e+03  &  2.0344e+05  & 1.5033e+06  \\
         $\kappa_2 (A_{11})$ &  1.6275e+05 &  8.8861e+06  &  2.6489e+10 &  1.4462e+12  \\
          $\normtwo{A_{11}}$  &  4.0343e+02 &  2.9810e+03 &  1.6275e+05 &  1.2026e+06  \\
          $\normtwo{A_{11}^{-1}}$  & 4.0343e+02 &  2.9810e+03  &  1.6275e+05  &  1.2025e+06  \\
           $\normtwo{A_{11}^{-1} - S}$ &      1.9563e-10 &  1.7270e-06 & 3.2124e-01 &  1.8255e+02  \\
           $\normtwo{A_{11}^{-1} - S}/\normtwo{A}$ &  3.8831e-13 &  4.6353e-10 &   1.5790e-06  &  1.2144e-04 \\
           $\normtwo{A-L_1 L_1^T}/\normtwo{A}$ &  3.8826e-13 &   4.6353e-10 &    1.5790e-06 &  1.2144e-04 \\
         $\normtwo{A-L_2 L_2^T}/\normtwo{A}$ &  5.6416e-17  & 0.0000e+00  & 1.8322e-16  &  7.7442e-17   \\
            $\normtwo{A^T  J A - J}$ &  1.2576e-11  &  4.6566e-10 &   3.2754e-06 &  1.0986e-04  \\
               $\normtwo{L_1^T  J L_1 - J}$ & 6.8689e-13  & 4.1809e-11  & 1.1162e-07    &  3.5415e-05 \\
                  $\normtwo{L_2^T  J L_2 - J}$ & 1.4451e-12  &  3.3716e-10 &  6.6639e-07  & 8.1038e-05   \\
	   \hline
  \end{tabular}
\end{center}
\end{table}

\medskip
The results are contained in Table $1$. 
Notice that the matrix $A_{11}$ is ill-conditioned, and  we observe a significant discrepancy
between the factorization  errors for Algorithms $W_1$ and $W_2$, as predicted (see Remark \ref{RemarkB}).
In contrast  to Algorithm $W_2$, Algorithm $W_1$ produces unstable results. We observe that the loss of symplecticity of the computed matrices $L_1$ and $L_2$ by these algorithms is almost on the same level.
\end{example}

\medskip

\begin{example}\label{Example2}
In the second  experiment we use the same matrix $S$ and repeat the calculations for the inverse of $A$ from Example \ref{example1}.
\medskip
\begin{table}\label{tabelka2}
\caption{The results for Example \ref{Example2} with  $A=(S^T \,S)^{-1}$, where $S$ is defined in (\ref{S}).}
\medskip
\begin{center}
\begin{tabular}{lccccc}
	\hline
	$\theta$ &  $3$ & $4$ & $6$ & $7$ \\
	\hline
 $\kappa_2 (A)$ & 2.5380e+05  & 1.3881e+07  & 4.1389e+10  & 2.2605e+12  \\
 $\normtwo{A}$ &   5.0379e+02  & 3.7257e+03 &  2.0344e+05  & 1.5033e+06  \\
 $\kappa_2 (A_{11})$ &  5.0198e+00  &  5.0027e+00 &  5.0001e+00  & 4.9995e+00 \\
  $\normtwo{A_{11}}$ & 5.0379e+02  & 3.7257e+03 &  2.0344e+05 &  1.5031e+06  \\
 $\normtwo{A_{11}^{-1}}$ &    9.9641e-03 &  1.3428e-03 &  2.4577e-05 &  3.3261e-06  \\
  $\normtwo{A_{11}^{-1} - S}$ &    3.2894e-13  & 1.9382e-12 &  2.4867e-11 &  9.6611e-11  \\
  $\normtwo{A_{11}^{-1} - S}/\normtwo{A}$ &    6.5293e-16  & 5.2022e-16 &  1.2223e-16 &  6.4266e-17 \\
 $\normtwo{A-L_1 L_1^T}/\normtwo{A}$ &    6.8168e-16  & 5.2143e-16 &  1.1574e-16 &  9.3479e-17  \\
  $\normtwo{A-L_2 L_2^T}/\normtwo{A}$ &  1.3800e-16 & 3.0514e-17  & 7.1528e-17  & 9.3479e-17  \\
    $\normtwo{A^T  J A - J}$ &   1.6648e-11  & 1.8615e-09 &  3.1158e-06  & 2.2403e-04  \\
    $\normtwo{L_1^T  J L_1 - J}$ & 2.8422e-14   & 4.5475e-13  & 7.2760e-12   & 5.8208e-11  \\
    $\normtwo{L_2^T  J L_2 - J}$ & 1.1251e-11  & 4.5853e-10   &  1.1479e-06  &  3.9094e-05  \\
   	\hline
	\end{tabular}
\end{center}
\end{table}
The results are contained in Table $2$. 
Notice that  the condition number of $A$ remains the same in both Examples  \ref{example1} and \ref{Example2}.
However,  in this particular case, the matrix $A_{11}$ is  perfectly well-conditioned,
and $\normtwo{A_{11}^{-1}-S}/\normtwo{A}$ is the same order as the  machine precision order.
Therefore, Algorithm $W_1$ produces numerically stable results, like Algorithm $W_2$. However, here $\normtwo{L_1^T  J L_1 - J}$ is much smaller  than $\normtwo{L_2^T  J L_2 - J}$.
On the other hand, $\normtwo{L_2^T  J L_2 - J}$ is smaller than the loss of symplecticity $\normtwo{A^T  J A - J}$ of the matrix $A$.
\end{example}

\medskip
\noindent
The following lemma can  be useful in creating examples of symmetric positive definite symplectic matrices:

\begin{lemma} [{\cite[Theorem 5.2]{Dopico}}]
\label{lemma4}
	Every  symmetric positive definite symplectic matrix  $A \in \mathbb R^{2n \times 2n}$
	can be written as $A=PDP^T$, where
	\[
	P=\left(
	\begin{array}{cc}
	I  &  0 \\
	H &   I
	\end{array}
	\right),
     	\quad
	D=\left(
	\begin{array}{cc}
	G  &  0 \\
	0 &   G^{-1}
	\end{array}
	\right),
	\]
	where $G\in \mathbb R^{n \times n} $ is symmetric positive definite and $H \in \mathbb R^{n \times n}$ is symmetric.
Notice that the matrix $A\in \mathbb R^{2n \times 2n}$ then equals
\[
A=\left(
	\begin{array}{cc}
	A_{11}  &  A_{12} \\
	A_{12}^T &   A_{22}
	\end{array}
	\right)=
    \left(
	\begin{array}{cc}
	G  &  GH \\
	HG &   HGH+G^{-1}
	\end{array}
	\right).
\]
\end{lemma}

\medskip

\begin{example}\label{example3}
In the third  experiment we apply  Lemma \ref{lemma4} for  $A\in \mathbb R^{2n \times 2n}$  having  all integer entries. Then the matrix $A$ is  symplectic in floating-point arithmetic.
We take $H=G^{-1}$, where $G \in \mathbb R^{n \times n}$ is symmetric positive definite. Then the matrix $A$ has the following form
\begin{equation}\label{new_A}
A=\left(
	\begin{array}{cc}
	G  &  I \\
	I &   2 G^{-1}
	\end{array}
	\right).
\end{equation}

Here the matrix $G$ is generated  from the Pascal matrix  by the following \textsl{MATLAB} code:
\begin{verbatim}
I=eye(n); Z=zeros(n);
G=pascal(n);
G=G(n:-1:1, n:-1:1);
inv_G=round(inv(G));
A=[G,I;I,2*inv_G];
\end{verbatim}

\medskip
\begin{table}\label{tabelka3}
\caption{The results for Example \ref{example3} for   $A\in \mathbb R^{2n \times 2n}$.}
\medskip
\begin{center}
\begin{tabular}{lccccc}
	\hline
	$n$ &  $6$ & $8$ & $10$ & $12$ \\
	\hline
 $\kappa_2 (A)$ & 4.4315e+05  & 8.2581e+07  & 1.6621e+10 & 3.5056e+12 \\
 $\normtwo{A}$ &  6.6569e+02  & 9.0874e+03  & 1.2892e+05 & 1.8723e+06  \\
 $\kappa_2 (A_{11})$ & 1.1079e+05   & 2.0645e+07   & 4.1552e+09  & 8.7639e+11  \\
  $\normtwo{A_{11}}$ &  3.3285e+02  & 4.5437e+03  & 6.4461e+04  & 9.3616e+05 \\
 $\normtwo{A_{11}^{-1}}$ &  3.3285e+02 &  4.5437e+03  &  6.4461e+04 & 9.3616e+05  \\
  $\normtwo{A_{11}^{-1} - S}$ &  4.6794e-11  & 6.6398e-08  & 2.4960e-05  & 3.9030e-02 \\
  $\normtwo{A_{11}^{-1} - S}/\normtwo{A}$ &  7.0293e-14   & 7.3066e-12  & 1.9361e-10  & 2.0846e-08  \\
 $\normtwo{A-L_1 L_1^T}/\normtwo{A}$ & 7.0331e-14  & 7.3066e-12 & 1.9361e-10  & 2.0846e-08  \\
  $\normtwo{A-L_2 L_2^T}/\normtwo{A}$ &  6.2068e-17  & 6.6499e-17  & 3.7068e-17  & 6.3200e-17  \\
    $\normtwo{A^T  J A - J}$ &  0  & 0 &  0  & 0 \\
    $\normtwo{L_1^T  J L_1 - J}$ & 1.2726e-14  & 6.0685e-13  & 1.7901e-12  &  5.4190e-11 \\
    $\normtwo{L_2^T  J L_2 - J}$ &   2.8402e-13  & 1.1703e-11  & 2.4951e-10  & 6.8509e-08  \\
   	\hline
	\end{tabular}
\end{center}
\end{table}

The results are contained in Table $3$. Note that although the matrix $A$ is {\bf pure symplectic} (i.e., $\normtwo{A^T  J A - J}=0$), Algorithm $W_1$ produces unstable solutions.

\end{example}

\medskip

\begin{example}\label{example7}
We applied  Lemma \ref{lemma4} for creating matrices $  A\in \mathbb R^{2n \times 2n} $  of the form $A=PDP^T$, where the matrices $G$ and $H$
are  generated  for  $n=1:100$  by the following \textsl{MATLAB} code:
\begin{verbatim}
rng(`default');
R = randn(n);
H = (R+R')/2;
G = R*R';
\end{verbatim}
Random matrices of entries are from the distribution $N(0,1)$. They were generated by \textsl{MATLAB} function "randn".
Before each usage the random number generator was reset to its initial state.

\medskip

\begin{figure}\label{obrazek1}
	\includegraphics[width=6cm]{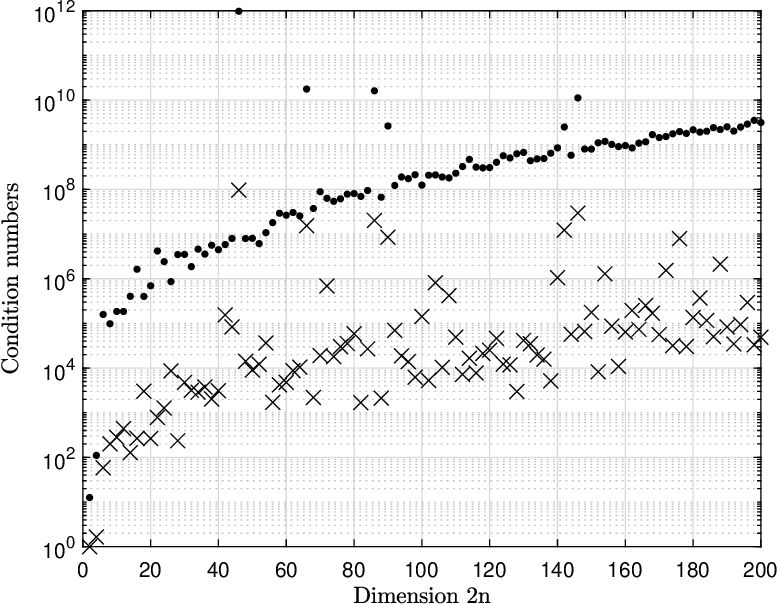}
	\includegraphics[width=6cm]{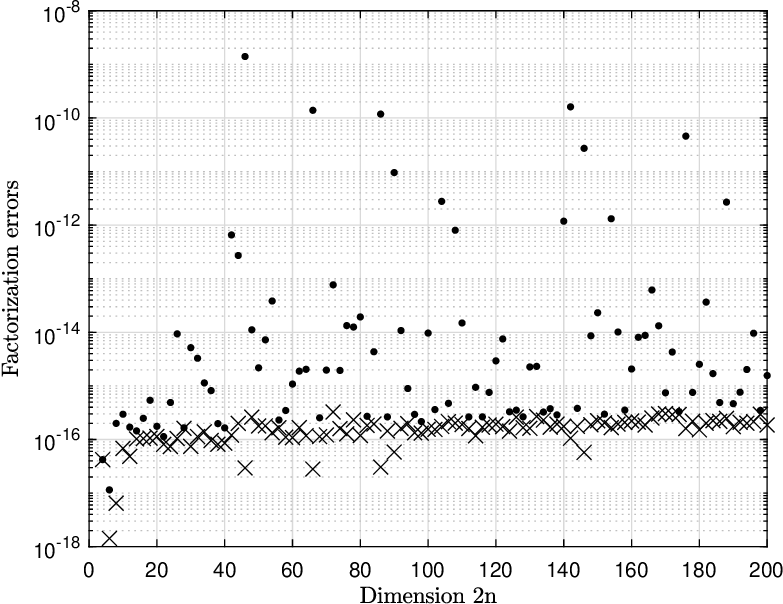}
	\caption{Condition numbers  $\kappa_2(A)$ (denoted by ".") and $\kappa_2(A_{11})$  (denoted "x")  for $A\in \mathbb R^{2n \times 2n}$ from  Example \ref{example7}.}
	\vskip-1.5ex
	\caption{Factorization error  $\normtwo{A-L_1 L_1^T}/\normtwo{A}$   for Algorithm $W_1$  (denoted by ".") and $\normtwo{A-L_2 L_2^T}/\normtwo{A}$ for Algorithm $W_2$  (denoted by "x ")  for  $A\in \mathbb R^{2n \times 2n}$ from  Example \ref{example7}.}
	\end{figure}

\vskip-1.5ex

\begin{figure}\label{obrazek3}
	\includegraphics[width=6cm]{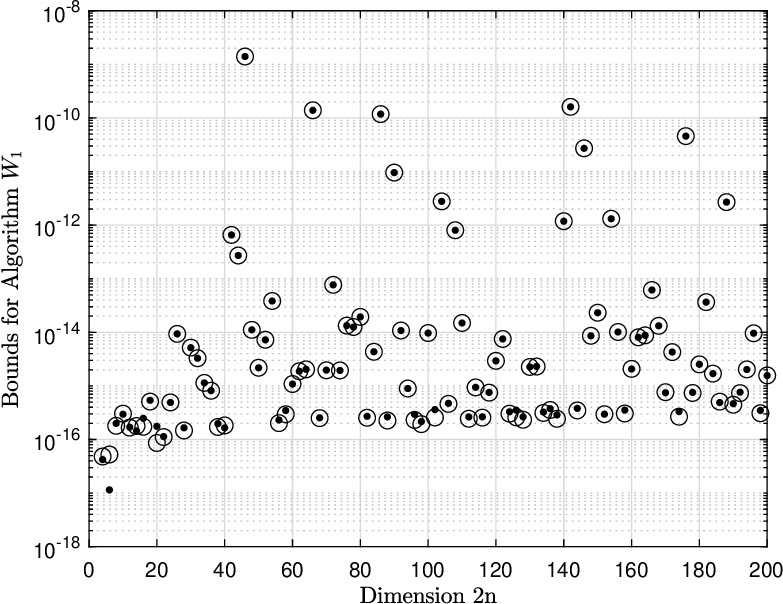}
	\includegraphics[width=6cm]{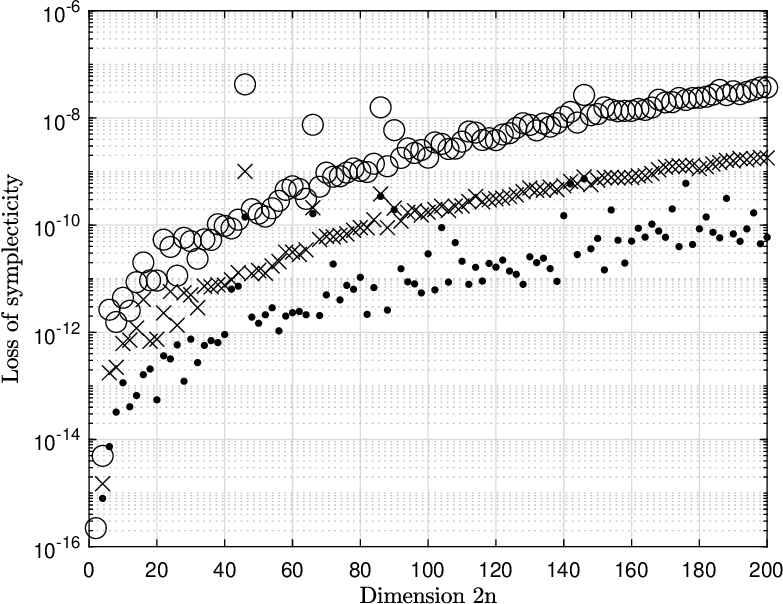}
	\caption{Factorization error   $\normtwo{A-L_1 L_1^T}/\normtwo{A}$   for Algorithm $W_1$  (denoted by ".") and $\normtwo{A_{11}^{-1}-S}/\normtwo{A}$  (denoted by "o") for  $A\in \mathbb R^{2n \times 2n}$ from  Example \ref{example7}.}
	\vskip-1.5ex
	\caption{Loss of  symplecticity $\normtwo{A^TJA-J}$   of  $A\in \mathbb R^{2n \times 2n}$ (denoted by "o") and $\normtwo{L_1^TJL_1-J}$  (denoted by ".") for  Algorithm $W_1$, and
	$\normtwo{L_2^TJL_2-J}$  (denoted by "x") for  Algorithm $W_2$  from  Example \ref{example7} .}
	\end{figure}

Figures \ref{obrazek1} and \ref{obrazek3}  report the values of the condition numbers  $\kappa_2(A)$ and  $\kappa_2(A_{11})$, the factorization error $\normtwo{A-L_1  L_1^T}/\normtwo{A}$ for Algorithm~$W_1$
 and the factorization error $\normtwo{A-L_2  L_2^T}/\normtwo{A}$ for Algorithm~$W_2$, respectively. The reader should pay attention to differences between factorization  errors in favor of Algorithm $W_2$.
It is clearly visible in Figure \ref{obrazek3} that the factorization error for Algorithm $W_1$ depends mainly on the quantity $\normtwo{A_{11}^{-1}-S}/\normtwo{A}$, as is also indicated by the bound (\ref{FirstW1}).
\end{example}

\section{Conclusions}

We have  studied  the numerical behavior of  Algorithms $W_1$ and $W_2$  for computing the symplectic $LL^T$ factorization
 and  derived the bounds for their factorization error. Theoretical results developed in the manuscript provide
 further insight into symplectic $LL^T$ factorization.
Through our analysis and numerical experiments, we have demonstrated that  Algorithm $W_2$ is numerically stable.
Algorithm $W_1$ is  cheaper than Algorithm $W_2$, but  it  is  unstable especially when the symmetric positive definite matrix $A$ is not exactly symplectic and the matrix $A_{11}$ is ill-conditioned.
Then its distance from the symplecticity is the dominant quantity for the factorization error in Algorithm $W_1$ that tries to preserve the structural properties of the symplectic factorization
although the matrix is only near symplectic.  We have shown that
Algorithm $W_1$ is guaranteed to be stable only if the principal submatrix matrix  $A_{11}\in \mathbb R^{n \times n}$ of a symmetric positive definite matrix $A\in \mathbb R^{2n \times 2n}$ is well-conditioned,
 and the distance to a symplecticity of $A$, measured as the departure of $A_{11}^{-1}$ from its Schur complement $S$, i.e.,  $\normtwo{A_{11}^{-1}-S}/\normtwo{A}$, is of the  order  of the  machine precision $\eps$, see   Remark \ref{RemarkB}.
The instability of Algorithm $W_1$ can be quite severe. Numerical tests confirm that our factorization error bounds of Algorithm $W_1$ are almost attainable.
In all our tests, $\normtwo{L_1^T  J L_1 - J}$ is  smaller  than $\normtwo{L_2^T  J L_2 - J}$. However,
$\normtwo{L_2^T  J L_2 - J}$ is smaller than the loss of symplecticity $\normtwo{A^T  J A - J}$ of the matrix $A$.

\medskip
{\bf ACKNOWLEDGEMENTS}

\medskip

The second author acknowledges support by Czech Academy of Sciences (RVO 67985840) and by the Grant Agency of the Czech Republic, Grant No. 23-06159S.
The second author is a
member of Ne\v{c}as Center for mathematical modelling.
The third author  acknowledges support from the EPSRC programme grant EP/R034826/1 and from the EPSRC research grant EP/V008129/1.

\end{document}